\documentclass[12pt]{article}
\usepackage{amssymb,amsmath,amsthm, amsfonts}
\usepackage{graphicx}
\usepackage{epsfig}
\usepackage{tikz}

\def\dref#1{(\ref{#1})}

\theoremstyle{plain}
\newtheorem{theorem}{Theorem}[section]
\newtheorem{lemma}{Lemma}[section]
\newtheorem{proposition}{Proposition}[section]
\newtheorem{corollary}{Corollary}[section]

\theoremstyle{definition}

\newtheorem{remark}{Remark}[section]

\numberwithin{equation}{section}

\begin{document}
\allowdisplaybreaks
\title{\large\bf Global boundedness and decay property of a three-dimensional Keller--Segel--Stokes system modeling coral fertilization}

\author
{\rm Jing Li\\
\it\small  College of Science, Minzu University of China,
  Beijing, 100081, P.R. China\\[2mm]
  \rm Peter Y.~H.~Pang\\
\it\small Department of Mathematics, National University of Singapore,\\
\it\small 10 Lower Kent Ridge Road, Republic of Singapore 119076\\[2mm]
\rm Yifu Wang\thanks{ \it\small Corresponding author. Email: {\tt wangyifu@bit.edu.cn}}\\
\it\small School of Mathematics and Statistics,
 Beijing Institute of Technology,\\
\it\small Beijing, 100081, P.R. China}
\vspace{-6em}
\date{}
\maketitle
\begin{abstract}
This paper is concerned with the four-component Keller--Segel--Stokes system modelling the fertilization process of corals:
\begin{equation*}
\left\{
\begin{array}{ll}
\rho_t+u\cdot\nabla\rho=\Delta\rho-\nabla\cdot(\rho\mathcal{S}(x,\rho,c)\nabla c)-\rho m, & \quad (x,t)\in \Omega\times (0,T),
\\
m_t+u\cdot\nabla m=\Delta m-\rho m,  & \quad (x,t)\in \Omega\times (0,T),
\\
c_t+u\cdot\nabla c=\Delta c-c+m,  & \quad (x,t)\in \Omega\times (0,T),
\\
u_t=\Delta u-\nabla P+(\rho+m)\nabla\phi,\quad \nabla\cdot u=0,  & \quad (x,t)\in \Omega\times (0,T)
\end{array}\right.
\end{equation*}
subject to the boundary conditions $\nabla c\cdot \nu=\nabla m\cdot \nu=(\nabla\rho-\rho \mathcal{S}(x,\rho,c)\nabla c)\cdot \nu=0$ and $u=0$, and suitably regular initial data $(\rho_0(x),m_0(x), c_0(x),u_0(x))$, where $T\in (0,\infty]$, $\Omega\subset\mathbb R^3$ is a  bounded domain with smooth boundary  $\partial\Omega$. This system describes the spatio-temporal dynamics of the population densities of sperm $\rho$ and egg $m$ under a chemotactic process facilitated by a chemical signal released by the egg with concentration $c$ in a fluid-flow environment $u$ modeled by the incompressible Stokes equation.
In this model, the chemotactic sensitivity tensor
$\mathcal{S}\in C^2(\overline\Omega\times [0,\infty)^2)^{3\times 3}$ satisfies
$|\mathcal{S}(x,\rho,c)|\leq C_S(1+\rho)^{-\alpha}$ with some $C_S>0$ and $\alpha\geq 0$. We will show that for  $\alpha\geq \frac 13$,
the solutions to the system are  globally bounded and  decay to a spatially homogeneous equilibrium
exponentially as time goes to infinity. In addition, we will also show that, for any  $\alpha\geq 0$, a similar result is valid when the initial data satisfy a certain smallness condition.

\end{abstract}\vspace{-1em}



\section{Introduction}

Chemotaxis, the directed movement caused by the concentration of certain chemicals,
 is ubiquitous in biology and ecology, and has a significant effect on pattern formation in numerous biological contexts \cite{Hillen,Maini}.
   The first mathematically rigorous studies
of chemotaxis were carried out by Patlak \cite{Patlak} and Keller--Segel \cite{Keller}. The
latter work involves the derivation of a system of PDEs, now known as the Keller--Segel system, which, despite its simple structure, was proved to have a lasting impact as a theoretical framework describing the collective behavior of populations under the influence of a chemotactic signal produced
by the populations themselves \cite{Bellomo,Herrero,Winkler2,Winkler1}. In contract to this well-understood Keller--Segel system, there seem to be
few theoretical results on nontrivial behavior in situations where the signal is not produced by the population, such as in oxygenotaxis processes of swimming aerobic bacteria \cite{Tuval}, or where the  signal production occurs by indirect processes, such as in glycolysis reaction and tumor invasion \cite{Dillon,Painter,Chaplain1}.

 In this paper, we study a chemotaxis--fluid system modelling coral fertilization.
Specifically, we are concerned with a Keller--Segel--Stokes system
\newpage
\begin{equation}\label{1.1}
\left\{
\begin{array}{ll}
\rho_t+u\cdot\nabla\rho=\Delta\rho-\nabla\cdot(\rho \mathcal{S}(x,\rho,c)\nabla c)-\rho m, &\!\!   (x,t)\in \Omega\times (0,T),
\\
m_t+u\cdot\nabla m=\Delta m-\rho m,  &\!\!  (x,t)\in \Omega\times (0,T),
\\
c_t+u\cdot\nabla c=\Delta c-c+m,  & \!\!   (x,t)\in \Omega\times (0,T),
\\
u_t=\Delta u-\nabla P+(\rho+m)\nabla\phi,\quad \nabla\cdot u=0,  &\!\!  (x,t)\in \Omega\times (0,T),
\\
\!\!(\nabla\rho-\rho \mathcal{S}(x,\rho,c)\nabla c)\cdot \nu\!\!=\!\nabla m\cdot \nu=\nabla c\cdot \nu=0,  u=0, & \!\!  \!\!  (x,t)\in \partial\Omega\times (0,T),\\
\!\!\!  \rho(x,0)\!=\!\rho_0(x),m(x,0)\!=\!m_0(x),c(x,0)\!=\!c_0(x),u(x,0)\!=\!u_0(x),&\!\!  x\in\Omega,
\end{array}\right.
\end{equation}
where $T\in (0,\infty]$, $\Omega\subset\mathbb R^3$ is a  bounded domain with smooth boundary  $\partial\Omega$,
the chemotactic sensitivity tensor $\mathcal{S}(x,\rho,c)=(s_{ij}(x,\rho,c))\in  C^2(\overline\Omega\times [0,\infty)^2)$, $i,j\in \{1,2,3\}$,
and $\phi\in W^{2,\infty}(\Omega)$.

This PDE system describes the phenomenon of coral broadcast spawning \cite{Espejo2,Espejo1,Kiselev1,Kiselev2}, where the sperm $\rho$ chemotactically moves toward the higher concentration of the chemical $c$ released by the egg $m$, while the egg $m$ is merely affected by random diffusion, fluid transport and degradation upon contact with the sperm. Meanwhile, the fluid flow vector $u$, modeling the ambient ocean environment, satisfies
a Stokes equation, where $P=P(x,t)$ represents the associated pressure,
and the buoyancy effect of the sperm and egg on the velocity, mediated through a given gravitational potential $\phi$, is taken into account.
We note that the use of the Stokes equation instead of the Navier--Stokes equation is justified by the observation 
that the fluid flow is relatively slow compared with the movement of the sperm and egg. We further note that
the sensitivity tensor $\mathcal{S}(x,\rho,c)$ may take values that are matrices possibly containing nontrivial off-diagonal entries, which reflects that the chemotactic migration may
not necessarily be oriented along the gradient of the chemical signal, but may rather involve rotational flux components (see \cite{Xue1,Xue2} for the detailed model derivation).

A two-component variant of \eqref{1.1} has been used in the mathematical study of coral broadcast spawning. Indeed, in \cite{Kiselev1,Kiselev2},
Kiselev and Ryzhik investigated  
the important effect of chemotaxis on the  
coral fertilization process  
via the Keller--Segel type system of the form  
 \begin{equation}\label{1.2}
\left\{
\begin{array}{ll}
\rho_t+u\cdot\nabla\rho=\Delta\rho-\chi\nabla\cdot(\rho \nabla c)-\mu\rho ^q,
\\
0=\Delta c+\rho
\end{array}\right.
\end{equation}
with a given regular solenoidal fluid flow vector $u$.
This model 
implicitly assumes that the densities of sperm and egg gametes are identical, and that the P\'{e}clet number for
the chemical concentration $c$ is small which allows us to ignore the effects of convection on $c$.  
The authors showed that, for the Cauchy problem in $\mathbb{R}^2$, the 
  total mass $ \int_{\mathbb{R}^2} \rho(x,t)dx$ can
become arbitrarily small with increasing $\chi$ in the
case $q > 2$ of supercritical reaction, 
whereas in the critical  case $q = 2$, a weaker
but related effect within finite time intervals is observed. Recently, 
 Ahn et al. \cite{JAEWOOK} established the global well-posedness
of regular solutions for the variant model of \eqref{1.2} with $c_t+u\cdot\nabla c=\Delta c-c+\rho$ instead  of $0=\Delta c+\rho$.  
 They also proved that  $ \int_{\mathbb{R}^d} \rho(x,t)dx$  $(d=2,3)$ asymptotically approaches a strictly positive constant $C(\chi)$ which tends to $0$ as $\chi\rightarrow \infty$. 

In \cite{Espejo3}, Espejo and Suzuki studied
the three-component variant of \eqref{1.1} 
\begin{equation}\label{1.3}
\left\{
\begin{array}{ll}
\rho_t+u\cdot\nabla\rho=\Delta\rho-\chi\nabla\cdot(\rho \mathcal{S}(x,\rho,c)\nabla c)-\mu \rho^2,
\\
c_t+u\cdot\nabla c=\Delta c-c+\rho,
\\
u_t+ \kappa (u\cdot \nabla) u=\Delta u-\nabla P+\rho\nabla\phi,
\\
 \nabla\cdot u=0
\end{array}\right.
\end{equation}
in the modeling of broadcast spawning when the interaction of chemotactic  movement of the gametes and  the surrounding fluid is not negligible. Here the coefficient $\kappa\in \mathbb{R}$  is related to the strength of nonlinear convection. In particular, when the fluid flow is slow, we can use the Stokes instead of
the Navier--Stokes equation, i.e., assume $\kappa = 0$ (see \cite{Difrancesco,Lorz}).
It should be mentioned that the chemotaxis--fluid  model with $c_t+u\cdot\nabla c=\Delta c-c\rho$ replacing the second equation in \eqref{1.3} has also been used to
describe the behavior of bacteria of the species Bacillus subtilis suspended in sessile water drops \cite{Tuval}. From the viewpoint of mathematical analysis, this chemotaxis--fluid system compounds the known difficulties in the study of fluid dynamics with the typical intricacies in the study of chemotaxis systems.
It has also been observed that when
$\mathcal{S}=\mathcal{S}(x,\rho,c)$ is a tensor, the corresponding chemotaxis--fluid system loses some energy-like structure, which plays a key role in the analysis of the scalar-valued
case.
Despite these challenges, some comprehensive results on the global-boundedness
and large time behavior of solutions are available in the literature (see \cite{Cao1,Li,Liu,Tao1,Wang1,Winkler3,Winkler5,Winkler6,Winklerp} for example).  It has been shown that
when $\mathcal{S}=\mathcal{S}(x,\rho,c)$ is a tensor fulfilling
\begin{equation}
|\mathcal{S}(x,\rho,c)|\leq\frac{C_{\mathcal{S}}}{(1+\rho)^{\alpha}}\quad\textrm{for some}~\alpha>0~\textrm{and}~ C_{\mathcal{S}}>0,\label{1.4}
\end{equation}
 the three-dimensional  system \dref{1.3} with $\mu=0$, $\kappa=0$ admits globally bounded weak solutions for $\alpha>1/2$ \cite{Wang1}, which is slightly stronger than
 the corresponding subcritical assumption $\alpha>1/3$ for the fluid-free system. As for $\alpha\geq 0$, 
when the suitably regular initial data $(\rho_0,c_0,u_0)$ fulfill a smallness
condition, 
\dref{1.3} with $\mu=0$, $\kappa=1$ possesses a global classical solution which  decays  to $(\bar{\rho}_0,\bar{\rho}_0,0)$ exponentially with $\bar{\rho}_0=\frac1{|\Omega|}\int_{\Omega} \rho_0(x)dx$ \cite{Yu}.

Removing the presupposition that
the densities of the sperm and egg coincide at each point,
Espejo and Suzuki \cite{Espejo2} looked at a simplified version of
\eqref{1.1} in two dimensions, namely,
\begin{equation}\label{1.5}
\left\{
\begin{array}{ll}
\rho_t+u\cdot\nabla\rho=\Delta\rho-\chi\nabla\cdot(\rho \nabla c)-\rho m,
\\
m_t+u\cdot\nabla m=\Delta m-\rho m,
\\
0=\Delta c+k_0(m-\displaystyle\frac1{|\Omega|}\int_{\Omega} m dx) ~\hbox{with}~ \int_{\Omega} c dx=0,
\end{array}\right.
\end{equation}
and showed that $\int_{\Omega} \rho_0(x)dx\geq \int_{\Omega} m_0 (x)dx$ implies that $m(x,t)$ vanishes asymptotically, while
$
\int_{\Omega} \rho(x,t)dx\rightarrow  \frac1{|\Omega|}(\int_{\Omega} \rho_0(x)dx- \int_{\Omega} m_0 (x) dx)
$
as $t\rightarrow \infty$, provided that $\chi $ is small enough and $u$ is low.
In two dimensions, Espejo and Winkler \cite{Espejo1} have  recently considered  the Navier--Stokes version of \eqref{1.1}:{\setlength\abovedisplayskip{4pt}
\setlength\belowdisplayskip{4pt}
\begin{equation}\label{1.6}
\left\{
\begin{array}{ll}
\rho_t+u\cdot\nabla\rho=\Delta\rho-\nabla\cdot(\rho \nabla c)-\rho m, \\
m_t+u\cdot\nabla m=\Delta m-\rho m,  \\
c_t+u\cdot\nabla c=\Delta c-c+m, \\
u_t +\kappa (u\cdot \nabla)=\Delta u-\nabla P+(\rho+m)\nabla\phi,\quad \nabla\cdot u=0,
\end{array}\right.
\end{equation}}
and established the global existence of classical solutions to the associated initial-boundary value problem, which tend towards a
spatially homogeneous equilibrium in the large time limit.

Motivated  by the above works, we shall  consider
the properties of solutions to the system \eqref{1.1} in the three-dimensional setting. In particular, we  shall show that the corresponding  solutions
converge to a spatially homogeneous equilibrium exponentially as $t\rightarrow \infty$ as well.

Throughout the rest of the paper, we shall assume that
\begin{equation}\label{1.7}
\left\{
\begin{array}{ll}
\rho_0\in C^0(\overline{\Omega}),~\rho_0\geq0 ~\hbox{and}~ \rho_0\not\equiv0,
\\
m_0\in C^0(\overline{\Omega}),~m_0\geq0 ~\hbox{and}~ m_0\not\equiv0,
\\
c_0\in W^{1,\infty}(\Omega),~c_0\geq0 ~\hbox{and}~ c_0\not\equiv0,
\\
u_0\in D(A^{\beta}) ~\hbox{for all}~ \beta\in(\frac34,1),
\end{array}\right.
\end{equation}
where $A$ denotes the realization of the Stokes operator in $L^2(\Omega)$.
Under these assumptions, we shall first  establish the existence of global bounded classical solutions to \eqref{1.1}:\vspace{-1em}
\begin{theorem}
Suppose that \eqref{1.4}, \eqref{1.7} hold with $\alpha>\frac13$. Then the system \eqref{1.1} admits a global classical solution $(\rho,m,c,u,P)$, which is uniformly bounded in the sense that for any $\beta\in(\frac34,1)$, there exists $K>0$ such that for all $ t\in(0,\infty)$
{\setlength\abovedisplayskip{4pt}
\setlength\belowdisplayskip{4pt}\begin{align}\label{1.8}
\|\rho(\cdot,t)\|_{L^\infty(\Omega)}+\|m(\cdot,t)\|_{L^\infty(\Omega)}+\|c(\cdot,t)\|_{W^{1,\infty}(\Omega)}
+\|A^\beta u(\cdot,t)\|_{L^2(\Omega)}\leq K.
\end{align}}
\end{theorem}

Then, we establish the large time behavior of these solutions as follows:\vspace{-1em}
\begin{theorem}
Under the assumptions of Theorem 1.1, the solutions given by Theorem 1.1 satisfy
{\setlength\abovedisplayskip{4pt}
\setlength\belowdisplayskip{4pt}$$
\rho(\cdot,t)\to \rho_\infty,~
m(\cdot,t)\to m_\infty,~
c(\cdot,t)\to m_\infty,~
u(\cdot,t)\to 0~
\hbox{in}~ L^\infty(\Omega) ~\hbox{as}~ t\to\infty.
$$}
Furthermore, when $\int_\Omega \rho_0\neq\int_\Omega m_0$, there exist $K>0$ and $\delta>0$ such that
{\setlength\abovedisplayskip{4pt}
\setlength\belowdisplayskip{4pt}
\begin{align}
\|\rho(\cdot,t)-\rho_\infty\|_{L^2(\Omega)}&\leq Ke^{-\delta t},\label{1.9}\\
\|m(\cdot,t)-m_\infty\|_{L^\infty(\Omega)}&\leq Ke^{-\delta t},\label{1.10}\\
\|c(\cdot,t)-m_\infty\|_{L^\infty(\Omega)}&\leq Ke^{-\delta t},\label{1.11}\\
\|u(\cdot,t)\|_{L^\infty(\Omega)}&\leq Ke^{-\delta t}, \label{1.12}
\end{align}}
where $\rho_\infty=\frac{1}{|\Omega|}\left\{\int_\Omega \rho_0-\int_\Omega m_0\right\}_+$, $m_\infty=\frac{1}{|\Omega|}\left\{\int_\Omega m_0-\int_\Omega \rho_0\right\}_+$.
\end{theorem}

According to the result for the related fluid--free system, the subcritical restriction  $\alpha>\frac13$ seems to be necessary
for the existence of global bounded solutions. However, for $\alpha\leq\frac13$, inspired by \cite{Cao1,Yu}, we investigate the existence of global bounded classical solutions and their large time behavior under a smallness assumption imposed on the
initial data, which can be stated as follows:\vspace{-1em}
\begin{theorem}
Suppose that \eqref{1.4} hold with $\alpha=0$ and $\int_{\Omega}\rho_0\neq\int_{\Omega}m_0$. Further, let $N=3$ and $p_0\in(\frac N2, \infty)$, $q_0\in(N,\infty)$ if $\int_{\Omega}\rho_0>\int_{\Omega}m_0$;
and $p_0\in(\frac {2N}3, \infty)$, $q_0\in(N,\infty)$ if $\int_{\Omega}\rho_0<\int_{\Omega}m_0$. Then there exists
$\varepsilon>0$ such that for any initial data $(\rho_0,m_0,c_0,u_0)$ fulfilling
\eqref{1.7} as well as {\setlength\abovedisplayskip{4pt}
\setlength\belowdisplayskip{4pt}
$$
\|\rho_0-\rho_\infty\|_{L^{p_0}(\Omega)}\leq\varepsilon,\quad  \|m_0-m_\infty\|_{L^{q_0}(\Omega)}\leq\varepsilon,
\quad\|\nabla c_0\|_{L^{N}(\Omega)}\leq\varepsilon, \quad\|u_0\|_{L^{N}(\Omega)}\leq\varepsilon,
$$
}
\eqref{1.1} possesses a global classical solution $(\rho,m,c,u,P)$.
Moreover, 
for any $\alpha_1$
$\in(0,\min\{\lambda_1, m_\infty+\rho_\infty\})$, $\alpha_2\in(0,\min\{\alpha_1,\lambda_1',1\})$, 
there exist constants $K_i$, $i=1,2,3,4$, such that for all $t\geq 1 $,   
{\setlength\abovedisplayskip{4pt}
\setlength\belowdisplayskip{4pt}
\begin{align*}
\|m(\cdot,t)-m_\infty\|_{L^\infty(\Omega)}\leq K_1e^{-\alpha_1 t},\quad
\|\rho(\cdot,t)-\rho_\infty\|_{L^\infty(\Omega)}\leq K_2e^{-\alpha_1 t},
\\
\|c(\cdot,t)-m_\infty\|_{W^{1,\infty}(\Omega)}\leq K_3e^{-\alpha_2t},
\quad
\|u(\cdot,t)\|_{L^\infty(\Omega)}\leq K_4 e^{-\alpha_2 t}.
\end{align*}}
Here $\lambda'_1$ is the first eigenvalue of $A$, and $\lambda_1$ is the first nonzero eigenvalue of $-\Delta$ on $\Omega$  under the Neumann boundary condition.
\end{theorem}
\vspace{-1.5em}
\begin{remark}
In Theorem 1.3, we have excluded the case $\int_{\Omega}\rho_0=\int_{\Omega}m_0$. Indeed, in this case, some results of Cao and Winkler \cite{Cao} suggest that exponential decay of solutions may not hold.
\end{remark}
\begin{remark} It is observed  that the similar result to Theorem 1.3 is also valid for the Navier--Stokes counterpart of (1.1) upon slight
modification of the definition of $T$ in \eqref{3.52} and \eqref{3.79}.
\end{remark}
\vspace{-1em}
As mentioned above, compared with the scalar sensitivity $\mathcal{S}$,
the system \eqref{1.1} with rotational tensor loses
 a favorable quasi-energy structure. For example, we note that the integral
 {\setlength\abovedisplayskip{4pt}
\setlength\belowdisplayskip{4pt}$$
 \int_\Omega \rho  ln  \rho+a\int_\Omega|\nabla c|^2+b \int_\Omega|u|^2
 $$}
with appropriate positive constants $a$ and $b$ plays a favorable entropy-like functional in deriving the bounds of solution to \eqref{1.6}. However, this will no longer be available in the present situation (see \cite{Espejo1}). To overcome this difficulty, our approach underlying the derivation of Theorem 1.1  will be based on 
the  estimate of the functional{\setlength\abovedisplayskip{4pt}
\setlength\belowdisplayskip{4pt}
\begin{align*}
\|\rho(\cdot,t)\|_{L^2(\Omega)}^2+\|u(\cdot,t)\|_{W^{1,2}(\Omega)}^2+\|\nabla c(\cdot,t)\|_{L^2(\Omega)}^2.
\end{align*}}
In addition, the proof of the exponential decay results in Theorem 1.2 relies on careful analysis of the functional
{\setlength\abovedisplayskip{4pt}
\setlength\belowdisplayskip{4pt}\begin{align*}
G(t):=\int_\Omega(\rho-\overline{\rho})^2+a\int_\Omega(m-\overline{m})^2+b\int_\Omega(c-\overline{c})^2
+c\int_\Omega\rho m
\end{align*}}
with suitable parameters $a,b,c>0$. Indeed, it can be seen that $G(t)$ satisfies the ODE:
$
G'(t)+\delta_1 G(t)\leq 0
$
for some $\delta_1>0$, and thereby the convergence rate of solutions in $L^2(\Omega)$ is established.
 At the same time, in comparison with the chemotaxis--fluid system  considered in \cite{Cao1,Yu}, due to
 {\setlength\abovedisplayskip{4pt}
\setlength\belowdisplayskip{4pt}
 $$\|e^{t\Delta}\omega\|_{L^p(\Omega)}\leq C_1\left(1+t^{-\frac N2(\frac1q-\frac1p)}\right)e^{-\lambda_1t}\|\omega\|_{L^q(\Omega)}
 $$}
 for all $\omega\in L^q(\Omega)$ with $\int_\Omega\omega=0 $, 
  $-\rho m$ in the first equation  of \eqref{1.1} gives rise to some difficulty in mathematical analysis despite its dissipative feature. Accordingly it requires a non-trivial application of the mass conservation of $\rho(x,t)-m(x,t)$.

 The plan of this paper is as follows:  In Section 2, we give
a  local  existence  result and some useful estimates. 
In Section 3,  in  the case of $\mathcal{S}$ vanishing on the boundary,  we investigate  the existence and large time behavior of global bounded classical solutions under the assumption of either $\alpha>\frac13$  or smallness of the initial data. In the last  section, on the basis of certain a priori estimates,  we give the proofs of our main results.  
\vspace{-2em}
\section{Preliminaries}  
\vspace{-1em}
In this section, we first recall a result on the local  existence of classical solutions, which
can be proved by a straightforward adaptation of well-known fixed point argument  
(see \cite{Winkler3} for example).\vspace{-1em}

\begin{lemma}\label{lemma2.5}
Suppose that \eqref{1.4}, \eqref{1.7} and
{\setlength\abovedisplayskip{4pt}
\setlength\belowdisplayskip{4pt}
\begin{equation}\label{2.1}
\mathcal{S}(x,\rho,c)=0, ~~(x,\rho,c)\in \partial\Omega\times [0,\infty)\times [0,\infty)
\end{equation}}
hold. Then there exist $T_{max}\in(0,\infty]$  and a classical solution
$(\rho,m,c,u,P)$ of \eqref{1.1} on $(0,T_{max})$. Moreover, $\rho,m,c$ are nonnegative in $\Omega\times(0,T_{max})$, and
if $T_{max}<\infty$, then for $\beta\in(\frac34,1)$,
$\lim_{t\to T_{max}}\left(\|\rho(\cdot,t)\|_{L^\infty(\Omega)}+\|m(\cdot,t)\|_{L^\infty(\Omega)}
+\|c(\cdot,t)\|_{W^{1,\infty}(\Omega)}+
\|A^{\beta}u(\cdot,t)\|_{L^2(\Omega)}\right)=\infty.$
This solution is unique, up to addition of constants to $P$.\vspace{-0.5em}
\end{lemma}

The following  elementary properties of the solutions in Lemma \ref{lemma2.5} are immediate consequences of the integration of
the first and second equations in \eqref{1.1}, as well as an application of the maximum principle to the second and third equations.\vspace{-1em}

\begin{lemma}\label{Lemma2.6}
Suppose that \eqref{1.4}, \eqref{1.7} and \eqref{2.1} hold. Then for all $t\in(0,T_{max})$, the solution of \eqref{1.1} from Lemma \ref{lemma2.5} satisfies
{\setlength\abovedisplayskip{4pt}
\setlength\belowdisplayskip{4pt}\begin{align}
&\|\rho(\cdot,t)\|_{L^1(\Omega)}\leq\|\rho_0\|_{L^1(\Omega)},\quad\|m(\cdot,t)\|_{L^1(\Omega)}\leq\|m_0\|_{L^1(\Omega)},
\label{2.2}
\\
&\int_0^t\|\rho(\cdot,s)m(\cdot,s)\|_{L^1(\Omega)}ds\leq\min\{\|\rho_0\|_{L^1(\Omega)},\|m_0\|_{L^1(\Omega)}\},
\label{2.3}
\\
&\|\rho(\cdot,t)\|_{L^1(\Omega)}-\|m(\cdot,t)\|_{L^1(\Omega)}=\|\rho_0\|_{L^1(\Omega)}-\|m_0\|_{L^1(\Omega)},
\label{2.4}
\\
&\|m(\cdot,t)\|_{L^2(\Omega)}^2+2\int_0^t\|\nabla m(\cdot,s)\|_{L^2(\Omega)}^2ds\leq \|m_0\|_{L^2(\Omega)}^2,
\label{2.5}
\\
&\|m(\cdot,t)\|_{L^\infty(\Omega)}\leq\|m_0\|_{L^\infty(\Omega)},
\label{2.6}
\\
&\|c(\cdot,t)\|_{L^\infty(\Omega)}\leq\max\{\|m_0\|_{L^\infty(\Omega)},\|c_0\|_{L^\infty(\Omega)}\}.
\label{2.7}
\end{align}}
\end{lemma}
\vspace{-2em}
\section {Proof of Theorems  for $\mathcal{S}=0$ on $\partial\Omega$} \vspace{-1em} 

In this section, we shall  consider the case in which  besides \eqref{1.4}, the sensitivity satisfies
 $\mathcal{S}=0$ on $\partial\Omega$. Under this hypothesis, the boundary condition for $\rho$ in \eqref{1.1} actually  reduces to the homogeneous
Neumann condition $\nabla \rho\cdot \nu=0$.\vspace{-1.5em}

\subsection{Global boundedness for $\mathcal{S}=0$ on $\partial\Omega$}
\vspace{-0.5em}

\begin{lemma}\label{lemma3.1}
Suppose that \eqref{1.4}, \eqref{1.7}, \eqref{2.1} hold with $\alpha>\frac13$. Then for any $\varepsilon>0$, there exists $K(\varepsilon)>0$ such that, for all $t\in(0,T_{max})$, the solution of \eqref{1.1} satisfies
\begin{align}\label{3.1}
&\frac{d}{dt}\|\rho(\cdot,t)\|_{L^2(\Omega)}^2+\frac12\|\nabla \rho(\cdot,t)\|_{L^2(\Omega)}^2\leq \varepsilon\|\Delta c(\cdot,t)\|_{L^2(\Omega)}^2+K(\varepsilon).
\end{align}

\end{lemma}

\proof Multiplying the first equation of \eqref{1.1} by $\rho$, we obtain
\begin{align}\label{3.2}
\frac12\frac{d}{dt}\!\int_\Omega\rho^2\!+\!\int_\Omega|\nabla\rho|^2\!=\!\!\int_\Omega\rho \mathcal{S}(x,\rho,c)\nabla\rho\nabla c\!\!-\!\int_\Omega\rho^2m
\leq\frac12\int_\Omega|\nabla\rho|^2+\frac{C_S^2}{2}\int_\Omega\frac{\rho^2}{(1\!+\!\rho)^{2\alpha}}|\nabla c|^2.
\end{align}
Now we estimate the term $\frac{C_S^2}{2}\int_\Omega\frac{\rho^2}{(1+\rho)^{2\alpha}}|\nabla c|^2$ in the right hand side of \eqref{3.2}.
In fact, if $\alpha\geq\frac34$,
\begin{align}\label{3.3}
\frac{C_S^2}{2}\int_\Omega\frac{\rho^2}{(1+\rho)^{2\alpha}}|\nabla c|^2\leq\varepsilon\int_\Omega|\nabla c|^4+K(\varepsilon),
\end{align}
while for $\alpha\in\left(\frac13,\frac34\right)$,
\begin{align}\label{3.4}
\frac{C_S^2}{2}\int_\Omega\frac{\rho^2}{(1+\rho)^{2\alpha}}|\nabla c|^2&\leq\frac{C_S^2}{2}\int_\Omega\rho^{2-2\alpha}|\nabla c|^2
\leq \frac{C_S^4}{16\varepsilon}\int_\Omega\rho^{4-4\alpha}+\varepsilon\int_\Omega|\nabla c|^4.
\end{align}
On the other hand, by Lemma \ref{Lemma2.6} and the Gagliardo--Nirenberg inequality, we get
{\setlength\abovedisplayskip{4pt}
\setlength\belowdisplayskip{4pt}\begin{align}\label{3.6}
\int_\Omega|\nabla c|^4& \leq C_{GN}\left\{\|\Delta c\|_{L^2(\Omega)}^2\|c\|_{L^\infty(\Omega)}^2+\|c\|_{L^\infty(\Omega)}^4\right\}
\leq C_{GN}'(\|\Delta c\|_{L^2(\Omega)}^2+1)
\end{align}}
and
{\setlength\abovedisplayskip{4pt}
\setlength\belowdisplayskip{4pt}
\begin{align*}
\int_\Omega|\rho|^{4-4\alpha}=\|\rho\|_{L^{4-4\alpha}(\Omega)}^{4-4\alpha}
&\leq C_{GN}\left\{\|\nabla\rho\|_{L^2(\Omega)}^{(4-4\alpha)\lambda_2}\|\rho\|_{L^1(\Omega)}^{(4-4\alpha)(1-\lambda_2)}
+\|\rho\|_{L^1(\Omega)}^{4-4\alpha}\right\}
\end{align*}}
with $\lambda_2=\frac{6(3-4\alpha)}{5(4-4\alpha)}$. Due to  $\alpha\in\left(\frac13,\frac34\right)$, we have $(4-4\alpha)\lambda_2<2$ and thus
{\setlength\abovedisplayskip{4pt}
\setlength\belowdisplayskip{4pt}
\begin{align}\label{3.7}
\frac{C_S^4}{16\varepsilon}\int_\Omega|\rho|^{4-4\alpha}\leq \frac14\int_\Omega|\nabla\rho|^2+K_1
\end{align}}
by the Young inequality. Combining \eqref{3.2}--\eqref{3.7}, we readily have \eqref{3.1}.
\begin{lemma} \label{lemma3.2}
Under the assumptions  of Lemma \ref{lemma3.1}, there exists a
positive constant $C = C(m_0,c_0)$ such that for all $t\in(0,T_{max})$, the solution of \eqref{1.1}  satisfies
\begin{align}
\frac{d}{dt}\|\nabla c(\cdot,t)\|_{L^2(\Omega)}^2+2\|\nabla c(\cdot,t)\|_{L^2(\Omega)}^2+\|\Delta c(\cdot,t)\|_{L^2(\Omega)}^2\leq K(\|\nabla u(\cdot,t)\|_{L^2(\Omega)}^2+1).\label{3.8}
\end{align}

\end{lemma}

\proof Multiplying the $c$-equation of \eqref{1.1} by $-\Delta c$, we obtain
\begin{align}
\frac12\frac{d}{dt}\int_\Omega|\nabla c|^2+\int_\Omega|\Delta c|^2+\int_\Omega|\nabla c|^2
&\leq-\int_\Omega m\Delta c+\int_\Omega (u\cdot\nabla c)\Delta c\label{3.9}
\\
&\leq\int_\Omega|m|^2+\frac14\int_\Omega|\Delta c|^2-\int_\Omega\nabla c\cdot(\nabla u\cdot\nabla c)\nonumber
\\
&\leq\|m\|_{L^2(\Omega)}^2+\frac14\|\Delta c\|_{L^2(\Omega)}^2+\frac{1}{2\varepsilon}\|\nabla u\|_{L^2(\Omega)}^2+
\frac{\varepsilon}{2}\|\nabla c\|_{L^4(\Omega)}^4.\nonumber
\end{align}
By \eqref{3.6} and taking $\varepsilon=\frac{1}{2C_{GN}'}$ in the above  inequality, we have
\begin{align*}
\frac12\frac{d}{dt}\int_\Omega|\nabla c|^2+\frac12\int_\Omega|\Delta c|^2+\int_\Omega|\nabla c|^2\leq\|m\|_{L^2(\Omega)}^2+C_{GN}'\|\nabla u\|_{L^2(\Omega)}^2+\frac14,
\end{align*}
which along with \eqref{2.5} readily ensures the validity of \eqref{3.8}.

\begin{lemma}
Under the assumptions  of Lemma \ref{lemma3.1}, the solution of \eqref{1.1} satisfies
\begin{align}
\frac{d}{dt}\|u(\cdot,t)\|_{L^2(\Omega)}^2+\|\nabla u(\cdot,t)\|_{L^2(\Omega)}^2\leq & K\left(\|\rho(\cdot,t)\|_{L^2(\Omega)}^2+1\right),\label{3.10}
\\
\frac{d}{dt}\|\nabla u(\cdot,t)\|_{L^2(\Omega)}^2+\|A u(\cdot,t)\|_{L^2(\Omega)}^2\leq & K\left(\|\rho(\cdot,t)\|_{L^2(\Omega)}^2+1\right)\label{3.11}
\end{align}
for all $t\in(0,T_{max})$ for a positive constant $K$.
\end{lemma}

\proof  Testing the $u$-equation in \eqref{1.1} by $u$,  using the   H\"{o}lder inequality and Poincar\'{e} inequality, we can get
{\setlength\abovedisplayskip{4pt}
\setlength\belowdisplayskip{4pt}
\begin{align*}
\frac12\frac{d}{dt}\int_\Omega|u|^2+\int_\Omega|\nabla u|^2&=\int_\Omega (\rho+m)\nabla\phi \cdot u
\\
&\leq \|\nabla\phi \|_{L^\infty (\Omega)}\|\rho+m\|_{L^2(\Omega)}\|u\|_{L^2(\Omega)}
\\
&\leq \frac 12\|\nabla u\|_{L^2(\Omega)}^2+K_1(\|\rho\|_{L^2(\Omega)}^2+\|m\|_{L^2(\Omega)}^2),
\end{align*}}
 which together with \eqref{2.5}  yields \eqref{3.10}.  Applying the Helmholtz projection $\mathcal{P}$ to the fourth equation in \eqref{1.1}, testing the resulting identity by $A u$ and using the Young inequality, we have
\begin{align*}
\frac12\frac{d}{dt}\int_{\Omega}|A^{\frac 12} u|^2+\int_\Omega|A u|^2&=-\int_\Omega \mathcal{P}[(\rho+m)\nabla\phi] \cdot A u
\\
&\leq \frac12\int_\Omega|A u|^2+K_2(\int_\Omega\rho^2+\int_\Omega m^2),
\end{align*}
which yields \eqref{3.11}, due to \eqref{2.5} and the fact that $\int_{\Omega}|\nabla u|^2=\int_\Omega|A^{\frac 12} u|^2 $.

\begin{lemma}\label{lemma 3.4} Under the assumptions  of Lemma \ref{lemma3.1}, one can find $C>0$ such that
for all  $t\in (0,T_{max})$, the solution of \eqref{1.1}  satisfies{\setlength\abovedisplayskip{4pt}
\setlength\belowdisplayskip{4pt}
\begin{align*}
\|\rho(\cdot,t)\|_{L^2(\Omega)}^2+\|\nabla c(\cdot,t)\|_{L^2(\Omega)}^2+\| u(\cdot,t)\|_{W^{1,2}(\Omega)}^2\leq K.
\end{align*}}
\end{lemma}
\proof  By the Gagliardo--Nirenberg inequality
\begin{align*}
\|\rho\|_{L^2(\Omega)}\leq C_{GN}\left(\|\nabla \rho\|_{L^2(\Omega)}^{\frac35}\|\rho\|_{L^1(\Omega)}^{\frac25}+\|\rho\|_{L^1(\Omega)}\right)
\end{align*}
and \eqref{3.1}, for any $\varepsilon>0$, there exists $K(\varepsilon)>0$ such that
\begin{align}\label{3.12}
&\frac{d}{dt}\|\rho\|_{L^2(\Omega)}^2+\|\rho\|_{L^2(\Omega)}^2+\frac14\|\nabla \rho\|_{L^2(\Omega)}^2\leq \varepsilon\|\Delta c\|_{L^2(\Omega)}^2+K_1(\varepsilon).
\end{align}
Adding \eqref{3.10}  and \eqref{3.11},  and by the Poincar\'{e} inequality, one can find constants  $K_i>0$, $i=2,3,4$, such that
\begin{equation}\label{3.13}
\begin{array}{rl}
\displaystyle\frac{d}{dt}(\|u\|_{L^2(\Omega)}^2+\|\nabla u\|_{L^2(\Omega)}^2)+K_2(\|u\|_{L^2(\Omega)}^2+\|\nabla u\|_{L^2(\Omega)}^2)
\leq & K_3\left(\|\rho\|_{L^2(\Omega)}^2+1\right)\\
\leq & \displaystyle\frac18\|\nabla\rho\|_{L^2(\Omega)}^2+K_4.
\end{array}
\end{equation}
Recalling \eqref{3.8}, we get
\begin{align}\label{3.14}
\frac{d}{dt}\|\nabla c\|_{L^2(\Omega)}^2+2\|\nabla c\|_{L^2(\Omega)}^2+\|\Delta c\|_{L^2(\Omega)}^2\leq K_5\left(\|\nabla u\|_{L^2(\Omega)}^2+1\right).
\end{align}
Now combining the above inequalities and choosing $\varepsilon=\frac{K_2}{2K_5}$, one can see that there exists some constant $K_6>0$ such that
{\setlength\abovedisplayskip{4pt}
\setlength\belowdisplayskip{4pt}\begin{align*}
Y(t):=\|\rho(\cdot,t)\|_{L^2(\Omega)}^2+\|u(\cdot,t)\|_{W^{1,2}(\Omega)}^2
+\frac{1}{\varepsilon}\|\nabla c(\cdot,t)\|_{L^2(\Omega)}^2
\end{align*}}
satisfies
$
Y'(t)+\delta Y(t)\leq K_6,
$
where $\delta=\min\{1,\frac{K_2}{2}\}$. Hence by an ODE comparison argument, we obtain
$Y(t)\leq K_7$ for some constant $K_7>0$ and thereby complete the proof.

With all of the above estimates at hand, we can now establish the global existence result in the
case  $\mathcal{S}=0$ on $\partial\Omega$.

{\bf Proof of Theorem 1.1 in the case $\mathcal{S}=0$ on $\partial\Omega$.}   To establish the existence of globally bounded classical solution, by the extensibility criterion in Lemma
\ref{lemma2.5}, we only need to show that
{\setlength\abovedisplayskip{6pt}
\setlength\belowdisplayskip{4pt}
\begin{align}\label{3.16}
\|\rho(\cdot,t)\|_{L^\infty(\Omega)}+\|m(\cdot,t)\|_{L^\infty(\Omega)}+\|c(\cdot,t)\|_{W^{1,\infty}(\Omega)}
+\|A^\beta u(\cdot,t)\|_{L^2(\Omega)}\leq K_1
\end{align}}
for all $t\in(0,T_{max})$ with some positive constant $K_1$ independent of $T_{max}$. To this end, by the estimate of Stokes operator (
Corollary 3.4 of \cite{WinklerCV}), we first get
{\setlength\abovedisplayskip{4pt}
\setlength\belowdisplayskip{4pt}
\begin{align}\label{3.17}
\|u\|_{L^\infty(\Omega)}\leq K_2\|u\|_{W^{1,5}(\Omega)}\leq K_3
\end{align}}
with positive constant $K_3>0$ independent of $T_{max}$, due to $\|\rho\|_{L^2(\Omega)}\leq K_4$ and $\|m\|_{L^\infty(\Omega)}\leq K_4$
from Lemma \ref{lemma 3.4}  and
Lemma \ref{Lemma2.6}, respectively.

By  Lemma 2.1 of \cite{Ishida1}, Lemma \ref{lemma 3.4} and the Young inequality, we have
\begin{align*}
\sup_{t\in(0,T_{max})}\|\nabla c\|_{L^\infty(\Omega)}&\leq K_5(1+\sup_{t\in(0,T_{max})}\|m-u\cdot\nabla c\|_{L^4(\Omega)})
\\
&\leq K_5(1+\sup_{t\in(0,T_{max})}(\|m\|_{L^4(\Omega)}+\|u\|_{L^6(\Omega)}\|\nabla c\|_{L^{12}(\Omega)}))
\\
&\leq K_5(1+\sup_{t\in(0,T_{max})}(\|m\|_{L^4(\Omega)}+\|u\|_{L^6(\Omega)}\|\nabla c\|_{L^2(\Omega)}^{\frac16}\|\nabla c\|_{L^\infty(\Omega)}^{\frac56}))
\\
&\leq K_6(1+\sup_{t\in(0,T_{max})}\|\nabla c\|_{L^\infty(\Omega)}^{\frac56}),
\end{align*}
which implies that $
\displaystyle\sup_{t\in(0,T_{max})}\|\nabla c(\cdot,t)\|_{L^\infty(\Omega)}\leq K_7$. Along 
with \eqref{2.7} this implies
$\|c(\cdot,t)\|_{W^{1,\infty}(\Omega)}\leq K_8$.
Furthermore, applying the variation-of-constants formula to the $\rho-$equation in \dref{1.1}  and by Lemma \ref{lemma 3.4},  we get
\begin{align*}
\|\rho\|_{L^{\infty}(\Omega)}
\leq
&
\|e^{t\Delta}\rho_0\|_{L^{\infty}(\Omega)}
+\int^t_0\|e^{(t-s)\Delta}\nabla\cdot( \rho
\mathcal{S}
\nabla c+\rho u)\|_{L^{\infty}(\Omega)}ds\\
\leq&\|\rho_0\|_{L^{\infty}(\Omega)}
+
C_4\int^t_0 (1+ (t-s)^{-\frac 78}) e^{-\lambda_1(t-s)}
\|\rho
\mathcal{S}
\nabla c+\rho u)\|_{L^4(\Omega)}ds\\
\leq&\|\rho_0\|_{L^{\infty}(\Omega)}
+
K_9\int^t_0 (1+ (t-s)^{-\frac 78}) e^{-\lambda_1(t-s)}
\|\rho\|_{L^4(\Omega)}ds\\
\leq&
\|\rho_0\|_{L^{\infty}(\Omega)}
+
K_9\int^t_0 (1+ (t-s)^{-\frac 78}) e^{-\lambda_1(t-s)}
\|\rho\|^{\frac 12}_{L^\infty(\Omega)}\|\rho\|^{\frac 12}_{L^2(\Omega)}ds\\
\leq&
\|\rho_0\|_{L^{\infty}(\Omega)}
+K_{10}\sup_{s\in(0,T_{max})}\|\rho\|^{\frac 12}_{L^\infty(\Omega)}
\end{align*}
with $K_{10}=K_9\displaystyle\sup_{t\in(0,T_{max})}\|\rho\|^{\frac 12}_{L^2(\Omega)}  \int^\infty_0 (1+ s^{-\frac 78}) e^{-\lambda_1s}
ds$,
where we have used $\nabla \cdot u=0$. Taking supremum on the left side of the above inequality  over $(0,T_{max})$, we obtain
 {\setlength\abovedisplayskip{4pt}
\setlength\belowdisplayskip{4pt}
$$
 \sup_{t\in(0,T_{max})}\|\rho\|_{L^{\infty}(\Omega)}\leq \|\rho_0\|_{L^{\infty}(\Omega)}
+K_{10}\sup_{t\in(0,T_{max})}\|\rho\|^{\frac 12}_{L^\infty(\Omega)},
 $$}
 and thereby $\displaystyle\sup_{t\in(0,T_{max})}\|\rho\|_{L^{\infty}(\Omega)}\leq K_{11}$ by the Young inequality.  Finally, by a
 straightforward argument, one  can find $K_{12}>0$  such that  $\displaystyle\sup_{t\in(0,T_{max})}\|A^\beta u\|_{L^2(\Omega)}\leq K_{12}$.
The boundedness estimate \eqref{3.16} is now a direct
consequence of the above inequalities and this completes the proof.\vspace{-0.5em}

\subsection{Large time behavior for $\mathcal{S}=0$ on $\partial\Omega$}\vspace{-0.5em}

This section is devoted to showing the large time  behavior of global solutions to \eqref{1.1} obtained in the above subsection.
In order to  derive the convergence properties of solution
with respect to the norm in $L^2(\Omega)$, we shall  make use of  the following lemma. In the sequel, we denote  $\overline{f}=\frac{1}{|\Omega|}\int_\Omega f(x)dx $.
\begin{lemma}(Lemma 4.6 of \cite{Espejo1})\label{lemma3.5a}
Let $\lambda>0$, $C>0$, and suppose that $y\in C^1([0,\infty))$ and $h\in C^0([0,\infty))$ are nonnegative functions satisfying
$
y'(t)+\lambda y(t)\leq h(t)$ for some $\lambda>0$ and all $t>0$. Then if $\int_0^\infty h(s)ds\leq C$, we have $y(t)\to0$ as $t\to\infty$.
\end{lemma}

By means of the testing procedure and the Young inequality,  we have  
\begin{align}
\frac{d}{dt}\int_\Omega(\rho-\overline{\rho})^2&=2\int_\Omega(\rho-\overline{\rho})(\Delta\rho-\nabla(\rho\mathcal{ S}(x,\rho,c)\nabla c)-u\cdot\nabla\rho-\rho m+\overline{\rho m})\label{1.31}
\\
&=-2\int_\Omega|\nabla\rho|^2+2\int_\Omega\rho \mathcal{S}(x,\rho,c)\nabla c\cdot\nabla\rho-2\int_\Omega(\rho-\overline{\rho})(\rho m-\overline{\rho m})\nonumber
\\
&\leq-\int_\Omega|\nabla\rho|^2+K_1\int_\Omega|\nabla c|^2-2\int_\Omega(\rho-\overline{\rho})\rho m,\nonumber
\end{align}\vspace{-2.5em}
\begin{align}
\frac{d}{dt}\int_\Omega(m-\overline{m})^2&=2\int_\Omega(m-\overline{m})(\Delta m-u\cdot\nabla m-\rho m+\overline{\rho m})\label{1.32}
\\
&=2\int_\Omega m(\Delta m-u\cdot\nabla m)-2\int_\Omega(m-\overline{m})(\rho m-\overline{\rho m})\nonumber
\\
&\leq-2\int_\Omega|\nabla m|^2-2\int_\Omega(m-\overline{m})\rho m,\nonumber
\end{align}
\vspace{-2.5em}
\begin{align}
\frac{d}{dt}\int_\Omega(c-\overline{c})^2&=2\int_\Omega(c-\overline{c})(\Delta c-u\cdot\nabla c-(c-\overline{c})+(m-\overline{m}))\label{1.33}
\\
&=2\int_\Omega c(\Delta c-u\cdot\nabla c)-2\int_\Omega(c-\overline{c})^2+2\int_\Omega(c-\overline{c})(m-\overline{m})\nonumber
\\
&\leq-2\int_\Omega|\nabla c|^2-\int_\Omega(c-\overline{c})^2+\int_\Omega(m-\overline{m})^2,\nonumber
\end{align}
\vspace{-2.5em}
\begin{align}
\frac{d}{dt}\int_\Omega|u|^2&=-2\int_\Omega|\nabla u|^2+2\int_\Omega(\rho+m)\nabla\phi\cdot u-2\int_\Omega\nabla P\cdot u\label{1.34}
\\
&=-2\int_\Omega|\nabla u|^2+2\int_\Omega(\rho-\overline{\rho}+m-\overline{m})\nabla\phi\cdot u\nonumber
\\
&\leq -2\int_\Omega|\nabla u|^2+K_2\left(\int_\Omega|\rho-\overline{\rho}+m-\overline{m}|^2\right)^{\frac12}\left(\int_\Omega|u|^2\right)^{\frac12}\nonumber
\\
&\leq -\int_\Omega|\nabla u|^2+K_3\left(\int_\Omega|\rho-\overline{\rho}|^2+\int_\Omega|m-\overline{m}|^2\right)\nonumber,
\end{align}
where $\nabla\cdot u=0 $, $u\mid_{\partial \Omega}=0$ and the boundedness of $u,\nabla\phi$ and $\mathcal{S}$ are used.

\begin{lemma}\label{Lemma3.7}
Under the assumptions  of Lemma \ref{lemma3.1},
{\setlength\abovedisplayskip{4pt}
\setlength\belowdisplayskip{4pt}\begin{align*}
\|(\rho-\overline{\rho})(\cdot,t)\|_{L^\infty(\Omega)}\to0\quad\hbox{as}~t\to\infty,
\\
\|(m-\overline{m})(\cdot,t)\|_{L^\infty(\Omega)}\to0\quad\hbox{as}~t\to\infty,
\\
\|(c-\overline{c})(\cdot,t)\|_{L^\infty(\Omega)}\to0\quad\hbox{as}~t\to\infty,
\\
\|u(\cdot,t)\|_{L^\infty(\Omega)}\to0\quad\hbox{as}~t\to\infty.
\end{align*}}
\end{lemma}

\proof From \eqref{1.31}--\eqref{1.34}, it follows that
\begin{align}
\frac{d}{dt}\int_\Omega(\rho-\overline{\rho})^2\leq & -\int_\Omega|\nabla\rho|^2+K_1\int_\Omega|\nabla c|^2+2\overline{\rho}\int_\Omega\rho m,\label{3.22}
\\
\frac{d}{dt}\int_\Omega(m-\overline{m})^2\leq & -2\int_\Omega|\nabla m|^2+2\overline{m}\int_\Omega\rho m,\label{3.23}
\\
\frac{d}{dt}\int_\Omega(c-\overline{c})^2\leq & -2\int_\Omega|\nabla c|^2-\int_\Omega(c-\overline{c})^2+\int_\Omega(m-\overline{m})^2,\label{3.24}
\\
\frac{d}{dt}\int_\Omega|u|^2\leq & -\int_\Omega|\nabla u|^2+K_3\left(\int_\Omega|\rho-\overline{\rho}|^2+\int_\Omega|m-\overline{m}|^2\right).\label{3.25}
\end{align}

Since $\int_\Omega|m-\overline{m}|^2\leq C_p \|\nabla m\|_{L^2(\Omega)}^2$ and $\int_0^\infty\int_\Omega\rho m\leq K_4$ by \eqref{2.3}, an application of
Lemma \ref{lemma3.5a}  to \eqref{3.23} yields
{\setlength\abovedisplayskip{4pt}
\setlength\belowdisplayskip{4pt} \begin{align}
\|m(\cdot,t)-\overline{m}(t)\|_{L^2(\Omega)}\to 0\quad \hbox{as}~t\to\infty\label{3.27}.
\end{align}}
Since
{\setlength\abovedisplayskip{4pt}
\setlength\belowdisplayskip{4pt}
 \begin{align}
 \int_0^\infty\int_\Omega|(m-\overline{m})|^2ds\leq C_p\int_0^\infty \|\nabla m\|_{L^2(\Omega)}^2ds \leq K_5,\label{3.333}
 \end{align}}
the  application of Lemma \ref{lemma3.5a}   to \eqref{3.24}  also yields
{\setlength\abovedisplayskip{4pt}
\setlength\belowdisplayskip{4pt}
\begin{align}
\|c(\cdot,t)-\overline{c}(t)\|_{L^2(\Omega)}\to 0\quad \hbox{as}~t\to\infty\label{3.28}
\end{align}}
and{\setlength\abovedisplayskip{4pt}
\setlength\belowdisplayskip{4pt}
\begin{align}
\int_0^\infty\|\nabla c\|_{L^2(\Omega)}^2\leq\int_0^\infty\int_\Omega|m-\overline{m}|^2+\int_\Omega|c_0-\overline{c_0}|^2 \leq K_6.\label{3.29}
\end{align}}
Furthermore, by \eqref{3.29}, $\int_\Omega|\rho-\overline{\rho}|^2\leq C_p \|\nabla \rho\|_{L^2(\Omega)}^2$ and $\int_0^\infty\int_\Omega\rho m\leq K_4$,
Lemma \ref{lemma3.5a} implies that
{\setlength\abovedisplayskip{4pt}
\setlength\belowdisplayskip{4pt}\begin{align}
& \|\rho(\cdot,t)-\overline{\rho}(t)\|_{L^2(\Omega)}\to 0 \quad
\hbox{as}~t\to\infty,\label{3.30}\\
& \int_0^\infty\|\rho-\overline{\rho}\|_{L^2(\Omega)}^2\leq  C_p\int_0^\infty\|\nabla \rho\|_{L^2(\Omega)}^2\leq K_7.\label{3.31}
\end{align}}
 Hence from \eqref{3.333}, \eqref{3.31}, $\int_\Omega|u|^2\leq C_p \|\nabla u\|_{L^2(\Omega)}^2$ and  Lemma \ref{lemma3.5a}, it follows that
{\setlength\abovedisplayskip{4pt}
\setlength\belowdisplayskip{4pt}\begin{align}
\|u(\cdot,t)\|_{L^2(\Omega)}\to 0\quad \hbox{as}~t\to\infty \label{3.32}
\end{align}}
as well as
$\int_0^\infty\|\nabla u\|_{L^2(\Omega)}^2\leq K_8.$

Now we turn the above convergence  in $L^2(\Omega)$  into $L^\infty(\Omega)$ with the help of  the higher regularity of the solutions.
Indeed, similar to the proof of $\|c(\cdot,t)\|_{W^{1,\infty}(\Omega)}\leq K$ in  Theorem 1.1 in the case $\mathcal{S}=0$ on $\partial\Omega$,
 $
 \|m(\cdot,t)\|_{W^{1,\infty}(\Omega)}\leq K_{10}
 $
 can be proved since $\|\rho(\cdot,t)\|_{L^\infty(\Omega)}+\|m(\cdot,t)\|_{L^\infty(\Omega)}\leq K_9$ for all $t>0$ in \eqref{3.16}.
 Hence from \eqref{3.16},
  there exists a constant $K_{11}>0$ such that
$
\|m(\cdot,t)-\overline{m}(t)\|_{W^{1,\infty}(\Omega)}\leq K_{11},~\|c(\cdot,t)-\overline{c}(t)\|_{W^{1,\infty}(\Omega)}\leq K_{11},~\|u(\cdot,t)\|_{W^{1,5}(\Omega)}\leq K_{11}
$
for all $t>1$. Therefore by \eqref{3.27}, \eqref{3.28} and \eqref{3.32},  the application of the  interpolation  inequality   yields
{\setlength\abovedisplayskip{4pt}
\setlength\belowdisplayskip{4pt}
\begin{align*}
& \|m-\overline{m}\|_{L^\infty(\Omega)}\leq C\left(\|m-\overline{m}\|_{W^{1,\infty}(\Omega)}^{\frac{3}{5}}
\|m-\overline{m}\|_{L^2(\Omega)}^{\frac{2}{5}}+\|m-\overline{m}\|_{L^2(\Omega)}\right)\to0\quad\hbox{as}~t\to\infty,\\
& \|c(\cdot,t)-\overline{c}(t)\|_{L^\infty(\Omega)}\to0, \quad
\|u(\cdot,t)\|_{L^\infty(\Omega)}\to0\quad\hbox{as}~t\to\infty.
\end{align*}}
In addition, similar to Lemma 4.4 in \cite{Espejo1} or Lemma 5.2 in  \cite{Cao1}, there exist $\vartheta\in (0,1)$ and constant $K_{12}>0$ such that
$
\|\rho\|_{C^{\vartheta,\frac \vartheta 2}(\overline{\Omega}\times [t,t+1])}\leq K_{12}
$
for all $t>1$, which along with \eqref{3.30} implies that
$
\|\rho(\cdot,t)-\overline{\rho}(t)\|_{C_{loc}(\overline{\Omega})}
\to 0 \quad
\hbox{as}~t\to\infty
$
and then, by the finite covering theorem,
$
\|\rho(\cdot,t)-\overline{\rho}(t)\|_{L^\infty({\Omega})}
\to 0 \quad
\hbox{as}~t\to\infty.
$

By very similar argument as in Lemma 4.2 of \cite{Espejo1}, we have
\begin{lemma}\label{Lemma3.8}
Under the assumptions  of Lemma \ref{lemma3.1},
{\setlength\abovedisplayskip{4pt}
\setlength\belowdisplayskip{4pt}
\begin{align*}
\overline{\rho}(t)\to\rho_\infty, \quad
\overline{m}(t)\to m_\infty,\quad
\overline{c}(t)\to m_\infty \quad\hbox{as}~t\to\infty
\end{align*}}
 with  $\rho_\infty=\{\overline{\rho_0}-\overline{m_0}\}_+$ and $m_\infty=\{\overline{m_0}-\overline{\rho_0}\}_+$.
\end{lemma}

\proof From \eqref{2.3} and \eqref{2.5}, we have{\setlength\abovedisplayskip{4pt}
\setlength\belowdisplayskip{4pt}
\begin{align}
& \int_{t-1}^t\|\rho m\|_{L^1(\Omega)}\to0\quad\hbox{as}~t\to\infty,\label{3.33}\\
& \int_{t-1}^t\|\nabla m\|_{L^2(\Omega)}^2\to0\quad\hbox{as}~t\to\infty.\label{3.34}
\end{align}}
On the other hand,
\begin{align*}
\int_{t-1}^t\|\rho m\|_{L^1(\Omega)}&=\int_{t-1}^t\int_\Omega\rho(m-\overline{m})
+\int_{t-1}^t\int_\Omega\rho\overline{m}
\\
&\geq-\int_{t-1}^t\|\rho(\cdot,s)\|_{L^2(\Omega)}\|m-\overline{m}\|_{L^2(\Omega)}+|\Omega|\int_{t-1}^t\overline{\rho}\cdot\overline{m}
\\
&\geq -K\int_{t-1}^t\|\nabla m\|_{L^2(\Omega)}+|\Omega|\int_{t-1}^t\overline{\rho}\cdot\overline{m}
\\
&\geq -K\left(\int_{t-1}^t\|\nabla m\|_{L^2(\Omega)}^2\right)^{\frac12}+|\Omega|\int_{t-1}^t\overline{\rho}\cdot\overline{m}.
\end{align*}
Inserting \eqref{3.33} and  \eqref{3.34} into the above inequality, we obtain
{\setlength\abovedisplayskip{4pt}
\setlength\belowdisplayskip{4pt}\begin{align}\label{3.35}
\int_{t-1}^t\overline{\rho}\cdot\overline{m}\to0\quad\hbox{as}~t\to\infty.
\end{align}}
Now if  $\overline{\rho_0}-\overline{m_0}\geq0$, \eqref{2.4} warrants that
$\overline{\rho}-\overline{m}\geq0$, which along with \eqref{3.35} implies that
{\setlength\abovedisplayskip{4pt}
\setlength\belowdisplayskip{4pt}\begin{align}
\int_{t-1}^t\overline{m}^2(s)ds\to0\quad\hbox{as}~t\to\infty.\label{1.53}
\end{align}}
Noticing that $\overline{m}(s)\geq\overline{m}(t)~\hbox{for all}~t\geq s$, we  have
$
0\leq\overline{m}(t)^2\leq \int_{t-1}^t \overline{m}^2(s) ds \to 0\quad\hbox{as}~t\to\infty,
$ and thus
$\overline{\rho}\to\rho_\infty ~\hbox{as}~t\to\infty$
due to \eqref{2.4}.
 By very similar argument,
one can see that
$\overline{\rho}\to 0 ~\hbox{as}~ t\to\infty$
and
$\overline{m}\to m_\infty ~\hbox{as}~t\to\infty$ in the case of $\overline{\rho_0}-\overline{m_0}<0$.
Finally,  it is observed that
$c(\cdot,t)\to m_\infty ~\hbox{in}~ L^2(\Omega)\,\,\hbox{as}~ t\to\infty
$
is also valid (see Lemma 4.7 of \cite{Espejo1} for example) and thus $\overline{c}(t)\to m_\infty~\hbox{as}~t\to\infty$ by the H\"{o}lder inequality.

Combining Lemma \ref{Lemma3.7} with Lemma \ref{Lemma3.8}, we have\vspace{-0.5em}
\begin{lemma}\label{lemma3.9}
Under the assumptions  of Lemma \ref{lemma3.1}, we have{\setlength\abovedisplayskip{4pt}
\setlength\belowdisplayskip{4pt}
$$\rho(\cdot,t)\to \rho_\infty,\;m(\cdot,t)\to m_\infty,\;
c(\cdot,t)\to m_\infty,\;
u(\cdot,t)\to 0~~
\hbox{in}~ L^\infty(\Omega) \,\, \hbox{as}\,\,t\to\infty.
$$}
\end{lemma}
Now we proceed to estimate the decay rate of
$\|\rho(\cdot,t)-\rho_\infty\|_{L^\infty(\Omega)}$, $\|m(\cdot,t)-m_\infty\|_{L^\infty(\Omega)}$,
$\|c(\cdot,t)-c_\infty\|_{L^\infty(\Omega)}$, and $\|u(\cdot,t)\|_{L^\infty(\Omega)}$ when $\int_\Omega\rho_0\neq \int_\Omega m_0$. To this end, we first consider
its decay rate in $ L^2(\Omega)$ based on a differential inequality.\vspace{-0.5em}

\begin{lemma}\label{Lemma3.10}
Under the assumptions  of Lemma \ref{lemma3.1} and $\int_\Omega \rho_0\neq\int_\Omega m_0$, for any $\varepsilon>0$, there exist constants $K(\varepsilon)>0$ and $t_\varepsilon>0$ such that for $t>t_\varepsilon$,
{\setlength\abovedisplayskip{4pt}
\setlength\belowdisplayskip{4pt}
\begin{align}
|\overline{\rho}(t)-\rho_\infty|+|\overline{m}(t)-m_\infty|\leq & K(\varepsilon)
e^{-(\rho_\infty+m_\infty-\varepsilon)t},\\
|\overline{c}(t)-m_\infty|\leq & K(\varepsilon)  e^{-\min\{1,(\rho_\infty+m_\infty-\varepsilon)\}t}.
\end{align}}
\end{lemma}
\proof For the case $\int_\Omega\rho_0>\int_\Omega m_0$, we have $\rho_\infty>0$ and $m_\infty=0$. By Lemma \ref{lemma3.9},
there exists $t_\varepsilon>0$ such that $\rho(x,t)\geq \rho_\infty-\varepsilon$ for $t>t_\varepsilon$ and $x\in \Omega$,
and thereby
$
\frac{d}{dt}\int_\Omega m=-\int_\Omega\rho m\leq-(\rho_\infty-\varepsilon)\int_\Omega m
$
for $t>t_\varepsilon$,  which implies that
$\overline{m}(t)\leq \overline{m_0} e^{-(\rho_\infty-\varepsilon) (t-t_\varepsilon)}$
 for $t>t_\varepsilon$.
Moreover, due to $\overline{\rho}=\overline{m}+\rho_\infty$ by \eqref{2.4},
we have
$
|\overline{\rho}(t)-\rho_\infty|=\overline{m}(t)\leq \overline{m_0}e^{-(\rho_\infty-\varepsilon) (t-t_\varepsilon)}\quad \hbox{for}~t>t_\varepsilon.
$
As for the case $\int_\Omega\rho_0<\int_\Omega m_0$, similarly we can prove that
$
|\overline{m}(t)-m_\infty|= \overline{\rho}\leq \overline{\rho_0}e^{-(m_\infty-\varepsilon) (t-t_\varepsilon)}.
$
for $t>t_\varepsilon$.
Furthermore, by the third equation of \eqref{1.1}, we have
$
\frac{d}{dt}\int_\Omega (c-m_\infty)=\int_\Omega (m-m_\infty)-\int_\Omega (c-m_\infty),
$
and thereby $
|\overline{c}(t)-m_\infty|\leq K(\varepsilon)e^{-\min\{1,\rho_\infty+m_\infty-\varepsilon\}t}.
$

{\bf Proof of Theorem 1.2 in the case $\mathcal{S}=0$ on $\partial\Omega$.}~ By Lemma \ref{Lemma3.7} and Lemma \ref{lemma3.9}, we have
{\setlength\abovedisplayskip{4pt}
\setlength\belowdisplayskip{4pt}
\begin{align*}
\rho(\cdot,t)-\overline{\rho}(t)\to 0,  ~
m(\cdot,t)-\overline{m}(t)\to 0, ~
\rho(\cdot,t)\to \rho_\infty,  ~
m(\cdot,t)\to m_\infty \quad\hbox{in}~L^\infty(\Omega)~\hbox{as}~t\to\infty,
\end{align*}}
which implies that for any $\varepsilon\in (0,\frac {\rho_\infty+ m_\infty}2 )$, there exists $t_\varepsilon>0$ such that $|\rho(\cdot,t)-\overline{\rho}(t)|<\varepsilon$, $|m(\cdot,t)-\overline{m}(t)|<\varepsilon$, $\rho(\cdot,t)+m(\cdot,t)\geq \rho_\infty+ m_\infty-\varepsilon$ for all $t>t_\varepsilon$ and $x\in \Omega$. Hence from \eqref{1.31}--\eqref{1.34}, we have
\begin{align}
& \frac{d}{dt}\int_\Omega(\rho-\overline{\rho})^2+\int_\Omega|\nabla\rho|^2\leq K_1\int_\Omega|\nabla c|^2+2\varepsilon\int_\Omega\rho m,\label{3.40}\\
& \frac{d}{dt}\int_\Omega(m-\overline{m})^2+2\int_\Omega|\nabla m|^2\leq 2\varepsilon\int_\Omega\rho m, \label{3.41}\\
& \frac{d}{dt}\int_\Omega(c-\overline{c})^2+2\int_\Omega|\nabla c|^2+\int_\Omega(c-\overline{c})^2\leq\int_\Omega(m-\overline{m})^2,\label{3.42}\\
& \frac{d}{dt}\int_\Omega|u|^2+\int_\Omega |\nabla u|^2\leq K_2\left(\int_\Omega(\rho-\overline{\rho})^2+\int_\Omega(m-\overline{m})^2\right)\label{3.43}
\end{align}
 for $t>t_\varepsilon$, as well as
\begin{align}
\frac{d}{dt}\int_\Omega \rho m&=\int_\Omega\left[\rho(\Delta m-u\cdot\nabla m-\rho m)+m(\Delta\rho-\nabla(\rho S(x,\rho,c)\nabla c)-u\cdot\nabla\rho-\rho m)\right]\nonumber\\
&=-2\int_\Omega\nabla\rho\nabla m-\int_\Omega(\rho u\cdot\nabla m+m u\cdot\nabla\rho)+\int_\Omega \rho S(x,\rho,c)\nabla c \cdot\nabla m-\int_\Omega \rho m^2-\int_\Omega \rho^2 m\nonumber
\\
&\leq\int_\Omega|\nabla\rho|^2+2\int_\Omega|\nabla m|^2-\int_\Omega u\cdot\nabla (\rho m)+K_3\int_\Omega |\nabla c|^2-\int_\Omega \rho m(\rho+m)\nonumber
\\
&\leq \int_\Omega|\nabla\rho|^2+2\int_\Omega |\nabla m|^2+K_3\int_\Omega |\nabla c|^2-\frac12(\rho_\infty+m_\infty)\int_\Omega \rho m,
\end{align}
where $\nabla\cdot u=0 $, $u\mid_{\partial \Omega}=0$ and the boundedness of $\rho$  are used.

On the other hand, by Poincare's inequality, there exists $C_P>0$ such that
\begin{align*}
\int_\Omega|\nabla\rho|^2\geq C_{P} \int_\Omega(\rho-\overline{\rho})^2,\quad\int_\Omega|\nabla m|^2\geq C_{P}\int_\Omega(m-\overline{m})^2,
\\
\quad\int_\Omega|\nabla c|^2\geq C_{P}\int_\Omega(c-\overline{c})^2,\quad\int_\Omega|\nabla u|^2\geq C_{P}\int_\Omega(u-\overline{u})^2.
\end{align*}
Therefore combining the above inequalities, and taking $\varepsilon<\frac{C_Pa(\rho_\infty+m_\infty)}{8(K_1+C_P}$ with
$a=\min\{\frac12,\frac{K_1}{4C_P},\frac{K_1}{K_3}\}$,  the functional
$
G(t):=\int_\Omega(\rho-\overline{\rho})^2+\frac{K_1}{C_P}\int_\Omega(m-\overline{m})^2+K_1\int_\Omega(c-\overline{c})^2
+a\int_\Omega\rho m
$
satisfies the ordinary differential inequality
$
\frac{d}{dt} G(t)+ \delta_1 G(t)\leq 0
$
with $\delta_1=\min\{\frac{C_P}{2}, 1,\frac{\rho_\infty+m_\infty}{4}\}$,
 which implies that
{\setlength\abovedisplayshortskip{2pt}
\setlength\belowdisplayshortskip{2pt}
\begin{align}\label{3.45}
\|\rho(\cdot,t)-\overline{\rho}\|_{L^2(\Omega)}+
\|m(\cdot,t)-\overline{m}\|_{L^2(\Omega)}+ 
\|c(\cdot,t)-\overline{c}\|_{L^2(\Omega)}\leq Ce^{-\frac{\delta_1}{2}t}.
\end{align}}
Moreover, by \eqref{3.45} and  \eqref{3.43},
$
\|u(\cdot,t)\|_{L^2(\Omega)}\leq Ce^{-\delta_2 t}
$ for some $\delta_2>0$.
At this position, combining \eqref{3.45} with Lemma \ref{Lemma3.10}, we can find $\delta_3>0$ such that{\setlength\abovedisplayskip{4pt}
\setlength\belowdisplayskip{4pt}
\begin{align}\label{3.47}
\|\rho(\cdot,t)-\rho_\infty\|_{L^2(\Omega)}+
\|m(\cdot,t)-m_\infty\|_{L^2(\Omega)}+
\|c(\cdot,t)-m_\infty\|_{L^2(\Omega)}\leq Ce^{-\delta_3 t}.
\end{align}}
Hence as in the proof of Lemma \ref{Lemma3.7}, we can obtain the decay estimates (1.9)--(1.12) by
an application of the interpolation  inequality, and thus the proof is complete.\vspace{-1em}

\subsection{Exponential decay under smallness condition}
\vspace{-0.5em}
In this subsection, we give the proof of Theorem 1.3 under the assumption that $\mathcal{S}=0$ on $\partial\Omega$.
The proof is divided into two cases (Proposition 3.1 and Proposition 3.2).\vspace{-1em}

\subsubsection{The case $\int_{\Omega}\rho_0>\int_{\Omega}m_0$}\vspace{-0.5em}

In this subsection we consider the case $\int_{\Omega}\rho_0>\int_{\Omega}m_0$, i.e.,
$\rho_\infty>0$, $m_\infty=0$.  

\begin{proposition} \label{Pr.3.1}
Suppose that \eqref{1.4} hold with $\alpha=0$ and $\int_{\Omega}\rho_0>\int_{\Omega}m_0$. Let $N=3$, $p_0\in(\frac N2, N)$, $q_0\in(N,\frac{Np_0}{N-p_0})$.
There exists
$\varepsilon>0$ such that for any initial data $(\rho_0,m_0,c_0,u_0)$ fulfilling
\eqref{1.7} as well as
{\setlength\abovedisplayskip{4pt}
\setlength\belowdisplayskip{4pt}
\begin{align*}
\|\rho_0-\rho_\infty\|_{L^{p_0}(\Omega)}\leq\varepsilon,\quad  \|m_0\|_{L^{q_0}(\Omega)}\leq\varepsilon,
\quad\|\nabla c_0\|_{L^{N}(\Omega)}\leq\varepsilon, \quad\|u_0\|_{L^{N}(\Omega)}\leq\varepsilon,
\end{align*}}
\eqref{1.1} admits a global classical solution $(\rho,m,c,u,P)$.
In particular, 
for any $\alpha_1\in(0,\min\{\lambda_1,\rho_\infty\})$, $\alpha_2\in(0,\min\{\alpha_1,\lambda_1',1\})$, 
there exist constants $K_i$, $i=1,2,3,4$,  such that for all $t\geq 1 $
{\setlength\abovedisplayskip{4pt}
\setlength\belowdisplayskip{4pt}\begin{align}
\|m(\cdot,t)\|_{L^\infty(\Omega)}\leq & K_1e^{-\alpha_1 t},\label{3.48}
\\
\|\rho(\cdot,t)-\rho_\infty\|_{L^\infty(\Omega)}\leq & K_2e^{-\alpha_1 t},\label{3.49}
\\
\|c(\cdot,t)\|_{W^{1,\infty}(\Omega)}\leq & K_3e^{-\alpha_2t},\label{3.50}
\\
\|u(\cdot,t)\|_{L^\infty(\Omega)}\leq & K_4 e^{-\alpha_2 t}.\label{3.51}
\end{align}}
\end{proposition}
Proposition 3.1 is the consequence of the following lemmas. In the proof of these lemmas, the constants $C_i>0$, $i=1,\ldots,10$, refer to those in
Lemma 1.2 and Lemma 1.3 of \cite{Winkler7}, Lemma 2.3 of \cite{Cao1}, Theorem 1 and Theorem 2 of \cite{Fujiwara}, respectively.
We first collect some easily verifiable observations in the following lemma:\vspace{-0.5em}
\begin{lemma}\label{Lemma3.11a}
Under the assumptions of Proposition 3.1 and $\sigma=\int_0^\infty\left(1+s^{-\frac{N}{2p_0}}\right)e^{-\alpha_1s}ds$, there exist $M_1>0,M_2>0$ and $\varepsilon>0$ such that
\begin{align}
&C_3+2 C_2C_{10} e^{(1+C_1+C_1|\Omega|^{\frac1{p_0}-\frac1{q_0}})\sigma} \leq \displaystyle \frac {M_2}4,
\quad \quad M_1 \varepsilon<1,\label{33.1}
\\
&12 C_2 C_{10}
(C_6\!+\! 4C_6C_9 C_{10} \|\nabla\Phi\|_{L^\infty(\Omega)}(M_1\!+\!C_1\!+\!C_1|\Omega|^{\frac1{p_0}\!-\!\frac1{q_0}}\!+4e^{(1\!+\!C_1\!+
\!C_1|\Omega|^{\frac1{p_0}\!-\!\frac1{q_0}})\sigma}))
\varepsilon<1,\label{33.2}
\\
&C_4 C_{10} C_SM_2 (e^{(1+C_1+C_1|\Omega|^{\frac1{p_0}-\frac1{q_0}})\sigma}+\rho_\infty|\Omega|^{\frac{1}{q_0}})\leq \displaystyle \frac {M_1}8,\label{33.3}
\\
&3 C_{10} C_4C_S(M_1+C_1+C_1|\Omega|^{\frac1{p_0}-\frac1{q_0}})M_2\varepsilon\leq \displaystyle \frac {M_1}8,\label{33.4}
\\
&3 C_{10}C_4 C_6(M_1+C_1+C_1|\Omega|^{\frac1{p_0}-\frac1{q_0}})
(1+2C_9 C_{10} \|\nabla\Phi\|_{L^\infty(\Omega)}(M_1+C_1+C_1|\Omega|^{\frac1{p_0}-\frac1{q_0}}+\nonumber
\\
& 4 e^{(1+C_1+C_1|\Omega|^{\frac1{p_0}-\frac1{q_0}})\sigma}))\varepsilon
\leq \displaystyle\frac{M_1}{4}.\label{33.5}
\end{align}
\end{lemma}

Let
\begin{align}
\!\!T\!\triangleq\!\sup\!\left\{
\!\widetilde{T}\!\in\!(0,T_{max})\!\left|
\begin{array}{l}
\!\|(\rho\!-\!m)(\cdot,t)\!-\!e^{t\Delta}(\rho_0\!-\!m_0)\|_{L^{\theta}(\Omega)}\!\leq\! M_1\varepsilon(1\!+\!t^{-\frac{N}{2}
(\frac{1}{p_0}-\frac{1}{\theta})})e^{-\alpha_1t}\\
\hbox{for all} \,\,\theta\in[q_0,\infty], ~t\in[0,\widetilde{T});
\\
 \!\|\nabla c(\cdot,t)\|_{L^\infty(\Omega)}\leq M_2\varepsilon(1+t^{-\frac12})e^{-\alpha_1t}~\hbox{for all}~t\in[0,\widetilde{T}).
\end{array}
\!\!\right.
\right\}\label{3.52}
\end{align}
By (1.7) and Lemma \ref{lemma2.5}, $T>0$ is well-defined. We first show $T=T_{max}$. To this end, we will show that all of the estimates
mentioned in \eqref{3.52} is valid with even smaller coefficients on the right hand side.
The derivation of these estimates will mainly rely on $L^p-L^q$ estimates for the Neumann heat semigroup
and the fact that the classical solutions on $(0,T_{max})$ can be represented as\vspace{-1em}
\begin{align}
(\rho-m)(\cdot,t)=& e^{t\Delta}(\rho_0-m_0)-\!\int_0^te^{(t-s)\Delta}(\nabla\cdot(\rho \mathcal{S}(x,\rho,c)\nabla c)+ u\cdot\nabla(\rho-m))(\cdot,s)ds,\label{3.53}
\\
m(\cdot,t)=& e^{t\Delta}m_0-\int_0^te^{(t-s)\Delta}(\rho m-u\cdot\nabla m)(\cdot,s)ds,\label{3.54}
\\
c(\cdot,t)=& e^{t(\Delta-1)}c_0+\int_0^te^{(t-s)(\Delta-1)}(m-u\cdot\nabla c)(\cdot,s)ds,\label{3.55}
\\
u(\cdot,t)=& e^{-tA}u_0+\int_0^te^{-(t-s)A}\mathcal{P}((\rho+m)\nabla\Phi)(\cdot,s)ds \label{3.56}
\end{align}
for all $t\in(0,T_{max})$ as per the variation-of-constants formula.
\begin{lemma}\label{Lemma3.11}
Under the assumptions of Proposition 3.1, for all $t\in(0,T)$ and $\theta\in[q_0,\infty]$,
\begin{align*}
\|(\rho-m)(\cdot,t)-\rho_\infty\|_{L^\theta(\Omega)}\leq M_3\varepsilon (1+t^{-\frac{N}{2}(\frac{1}{p_0}-\frac{1}{\theta})})e^{-\alpha_1t}.
\end{align*}

\end{lemma}

\proof Since $e^{t\Delta}\rho_\infty=\rho_\infty$ and $\int_\Omega (\rho_0-m_0-\rho_\infty)=0$, the definition of $T$ and
Lemma 1.3 of \cite{Winkler7} 
show that
{\setlength\abovedisplayskip{5pt}
\setlength\belowdisplayskip{4pt}
\begin{align*}
&\|(\rho-m)(\cdot,t)-\rho_\infty\|_{L^\theta(\Omega)}
\\
\leq&\|(\rho-m)(\cdot,t)-e^{t\Delta}(\rho_0-m_0)\|_{L^\theta(\Omega)}
+\|e^{t\Delta}(\rho_0-m_0-\rho_\infty)\|_{L^\theta(\Omega)}
\\
\leq& M_1\varepsilon (1+t^{-\frac{N}{2}(\frac{1}{p_0}-\frac{1}{\theta})})e^{-\alpha_1t}
+C_1(1+t^{-\frac{N}{2}(\frac{1}{p_0}-\frac{1}{\theta})})(\|\rho_0-\rho_\infty\|_{L^{p_0}(\Omega)}
+\|m_0\|_{L^{p_0}(\Omega)})e^{-\lambda_1t}
\\
\leq&M_3\varepsilon(1+t^{-\frac{N}{2}(\frac{1}{p_0}-\frac{1}{\theta})})e^{-\alpha_1t}
\end{align*}}
for all $t\in(0,T)$ and $\theta\in[q_0,\infty]$, where $M_3=M_1+C_1+C_1|\Omega|^{\frac1{p_0}-\frac1{q_0}}$.
\begin{lemma}\label{Lemma3.12}
Under the assumptions of Proposition 3.1, for any $k>1$, 
{\setlength\abovedisplayskip{4pt}
\setlength\belowdisplayskip{4pt}
\begin{align}\label{3.57}
\|m(\cdot,t)\|_{L^k(\Omega)}\leq M_4 \|m_0\|_{L^k(\Omega)} e^{-\rho_\infty t}\quad\hbox{for all}~ t\in(0,T)
\end{align}}
with $\sigma=\int_0^\infty(1+s^{-\frac{N}{2p_0}})e^{-\alpha_1s}ds$ and $M_4=e^{M_3\sigma \varepsilon }$.
\end{lemma}\vspace{-1em}
\proof Multiplying the $m$-equation in \eqref{1.1} by $km^{k-1}$ and integrating the result over $\Omega$, we get
$
\frac{d}{dt}\int_\Omega m^k\leq-k\int_\Omega\rho m^k
$ on $(0,T)$.
Since $-\rho\leq|\rho-m-\rho_\infty|-m-\rho_\infty\leq-\rho_\infty+|\rho-m-\rho_\infty|$,
Lemma \ref{Lemma3.11} yields
{\setlength\abovedisplayskip{4pt}
\setlength\belowdisplayskip{4pt}\begin{align*}
\frac{d}{dt}\int_\Omega m^k&\leq-k\rho_\infty\int_\Omega m^k+k\int_\Omega m^k|\rho-m-\rho_\infty|
\\
&\leq-k\rho_\infty\int_\Omega m^k+k\|\rho-m-\rho_\infty\|_{L^\infty(\Omega)}\int_\Omega m^k
\\
&\leq-k\rho_\infty\int_\Omega m^k+kM_3\varepsilon\left(1+t^{-\frac{N}{2p_0}}\right)e^{-\alpha_1t}\int_\Omega m^k
\end{align*}}
and thus
$
\int_\Omega m^k \leq\int_\Omega m_0^k \exp\{-k\rho_\infty t+kM_3\varepsilon \int_0^t(1+s^{-\frac{N}{2p_0}})e^{-\alpha_1s}ds\}
\leq \|m_0\|_{L^k(\Omega)} ^ke^{k(M_3\sigma\varepsilon -\rho_\infty t)}$.
 The assertion \eqref{3.57} follows immediately.\vspace{-0.5em}

\begin{lemma}\label{Lemma 3.14}
Under the assumptions of Proposition 3.1, there exists $M_3>0$ such that
$
\|u(\cdot,t)\|_{L^{q_0}(\Omega)}\leq M_5\varepsilon\left(1+t^{-\frac12+\frac{N}{2q_0}}\right) e^{-\alpha_2 t}
$ for all $t\in(0,T)$.
\end{lemma}

\proof For any given $\alpha_2<\lambda_1'$, we fix $ \mu\in (\alpha_2, \lambda_1')$.  By \eqref{3.56},
 Lemma 2.3 of \cite{Cao1}, Theorem 1 and Theorem 2 of \cite{Fujiwara}, we obtain
{\setlength\abovedisplayskip{4pt}
\setlength\belowdisplayskip{4pt}
\begin{align}
&\|u(\cdot,t)\|_{L^{q_0}(\Omega)}\nonumber
\\
\leq& C_6t^{-\frac N2\left(\frac1N-\frac{1}{q_0}\right)}e^{-\mu t}\|u_0\|_{L^N(\Omega)}
+\int_0^t\|e^{-(t-s)A}\mathcal{P}((\rho+m)\nabla\Phi)(\cdot,s)\|_{L^{q_0}(\Omega)}ds\label{3.58}
\\
\leq& C_6t^{-\frac N2\left(\frac 1N-\frac{1}{q_0}\right)}e^{-\mu t}\|u_0\|_{L^N(\Omega)}
+C_6\int_0^te^{-\mu(t-s)}\|\mathcal{P}((\rho+m-\overline{\rho+m})\nabla\Phi)(\cdot,s)\|_{L^{q_0}(\Omega)}ds\nonumber
\\
\leq& C_6 t^{-\frac 12+\frac{N}{2q_0}}e^{-\mu t}\|u_0\|_{L^N(\Omega)}
+C_6C_9\|\nabla\Phi\|_{L^\infty(\Omega)}\int_0^te^{-\mu(t-s)}\|(\rho+m-\overline{\rho+m})(\cdot,s)\|_{L^{q_0}(\Omega)}ds,\nonumber
\end{align}}
where $\mathcal{P}(\overline{\rho+m}\nabla\Phi)=\overline{\rho+m} \mathcal{P}(\nabla\Phi)=0$ is used.
On the other hand, due to  $\alpha_1<\rho_\infty$,  Lemma  \ref{Lemma3.11} and Lemma \ref{Lemma3.12} show that
{\setlength\abovedisplayskip{4pt}
\setlength\belowdisplayskip{4pt}
\begin{align}
&\|(\rho+m-\overline{\rho+m})(\cdot,s)\|_{L^{q_0}(\Omega)}\nonumber
\\
=&\|(\rho-m-\overline{\rho-m})(\cdot,s)+2(m-\overline{m})(\cdot,s)\|_{L^{q_0}(\Omega)}\label{1.50}
\\
\leq&\|(\rho-m-\rho_\infty)(\cdot,s)\|_{L^{q_0}(\Omega)}+2\|(m-\overline{m})(\cdot,s)\|_{L^{q_0}(\Omega)}\nonumber
\\
\leq&M_5'
\varepsilon(1+s^{-\frac{N}{2}(\frac{1}{p_0}-\frac{1}{q_0})})e^{-\alpha_1s}
\nonumber
\end{align}}
with $M_5'=M_3+4e^{M_3 \sigma \varepsilon}$.  Combining \eqref{3.58} with \eqref{1.50} and applying Lemma 1.2 of \cite{Winkler7},
we have
{\setlength\abovedisplayskip{4pt}
\setlength\belowdisplayskip{4pt}
\begin{align*}
&\|u(\cdot,t)\|_{L^{q_0}(\Omega)}
\\
\leq& C_6t^{-\frac 12+\frac{N}{2q_0}}e^{-\mu t}\|u_0\|_{L^N(\Omega)}
+C_6C_9\|\nabla\Phi\|_{L^\infty(\Omega)} M_5'\varepsilon\int_0^t(1+s^{-\frac N2\left(\frac{1}{p_0}-\frac{1}{q_0}\right)})e^{-\alpha_1s}e^{-\mu(t-s)}ds
\\
\leq& C_6 t^{-\frac12+\frac{N}{2q_0}}e^{-\mu t}\|u_0\|_{L^N(\Omega)}
+C_6C_9 C_{10} \|\nabla\Phi\|_{L^\infty(\Omega)}M_5'\varepsilon(1+t^{\min\{0,1-\frac N2\left(\frac{1}{p_0}-\frac{1}{q_0}\right)\}})e^{-\alpha_2t}
\\
\leq& C_6 t^{-\frac12+\frac{N}{2q_0}}e^{-\mu t}\varepsilon
+2 C_6C_9 C_{10} \|\nabla\Phi\|_{L^\infty(\Omega)}M_5'\varepsilon e^{-\alpha_2t}
\\
\leq& M_5\varepsilon (1+t^{-\frac12+\frac{N}{2q_0}})e^{-\alpha_2t},
\end{align*}}
where $M_5=C_6+ 2 C_6C_9 C_{10} \|\nabla\Phi\|_{L^\infty(\Omega)}M_5' $ and
 $\frac{N}{2}\left(\frac{1}{p_0}-\frac{1}{q_0}\right)<1$ is used.

\begin{lemma}\label{lemma3.15}
Under the assumptions of Proposition 3.1, for all $t\in(0,T)$,
\begin{align*}
\|\nabla c(\cdot,t)\|_{L^{\infty}(\Omega)}\leq \frac{M_2}{2}\varepsilon \left(1+t^{-\frac12}\right) e^{-\alpha_1 t}.
\end{align*}

\end{lemma}

\proof By \eqref{3.55} and Lemma 1.3 of \cite{Winkler7}, we have
\begin{align}
\|\nabla c(\cdot,t)\|_{L^\infty(\Omega)}
&\leq\|e^{t(\Delta-1)}\nabla c_0\|_{L^\infty(\Omega)}+\int_0^t\|\nabla
e^{(t-s)(\Delta-1)}(m-u\cdot\nabla c)(\cdot,s)\|_{L^\infty(\Omega)}ds\nonumber
\\
&\leq C_3(1+t^{-\frac12})e^{-(\lambda_1+1)t}\|\nabla c_0\|_{L^N(\Omega)}+\int_0^t\|\nabla
e^{(t-s)(\Delta-1)}m(\cdot,s)\|_{L^\infty(\Omega)}ds\nonumber
\\
&\quad+\int_0^t\|\nabla
e^{(t-s)(\Delta-1)}u\cdot\nabla c(\cdot,s)\|_{L^\infty(\Omega)}ds.\label{3.60}
\end{align}
Now we estimate the last two integrals on the right hand side. From Lemma 1.3 of \cite{Winkler7},  Lemma \ref{Lemma3.12} with $k=q_0$, Lemma 1.2 of \cite{Winkler7} and the fact that $q_0>N$, it follows that{\setlength\abovedisplayskip{4pt}
\setlength\belowdisplayskip{4pt}
\begin{align}\label{3.61}
\int_0^t\|\nabla
e^{(t-s)(\Delta-1)}m\|_{L^\infty(\Omega)}ds
\leq&C_2\int_0^t(1+(t-s)^{-\frac12-\frac{N}{2q_0}})e^{-(\lambda_1+1)(t-s)}\|m\|_{L^{q_0}(\Omega)}ds
\\
\leq&C_2 M_4\varepsilon
\int_0^t(1+(t-s)^{-\frac12-\frac{N}{2q_0}})e^{-(\lambda_1+1)(t-s)}
e^{-\rho_\infty s}ds
\nonumber
\\
\leq& C_2C_{10} M_4
(1+t^{\min\{0,\frac12-\frac{N}{2q_0}\}})\varepsilon e^{-\alpha_1t}\nonumber
\\
\leq& 2 C_2C_{10} M_4(1+t^{-\frac12})\varepsilon e^{-\alpha_1 t}.\nonumber
\end{align}}
On the other hand, by  Lemma 2.3 of \cite{Cao1}, Lemma \ref{Lemma 3.14} and the definition of $T$, we obtain{\setlength\abovedisplayskip{4pt}
\setlength\belowdisplayskip{4pt}
\begin{align}
&\int_0^t\|\nabla
e^{(t-s)(\Delta-1)}u\cdot\nabla c\|_{L^\infty(\Omega)}ds\nonumber
\\
\leq&C_2\int_0^t(1+(t-s)^{-\frac12-\frac{N}{2q_0}})e^{-(\lambda_1+1)(t-s)}\|u\cdot\nabla c\|_{L^{q_0}(\Omega)}ds\label{3.62}
\\
\leq&C_2\int_0^t(1+(t-s)^{-\frac12-\frac{N}{2q_0}})e^{-(\lambda_1+1)(t-s)}\|u\|_{L^{q_0}(\Omega)}\|\nabla c\|_{L^\infty(\Omega)}ds\nonumber
\\
\leq& C_2 M_5M_2 \varepsilon^2 \int_0^t(1+(t-s)^{-\frac12-\frac{N}{2q_0}})e^{-(\lambda_1+1)(t-s)}(1+s^{-\frac12+\frac{N}{2q_0}}) (1+s^{-\frac12}) e^{-(\alpha_1+\alpha_2) s}ds \nonumber
\\
\leq&3C_2M_5M_2\varepsilon^{2}\int_0^te^{-(\lambda_1+1)(t-s)}
 e^{-(\alpha_1+\alpha_2) s}(1+(t-s)^{-\frac12-\frac{N}{2q_0}})(1+s^{-1+\frac{N}{2q_0}})ds\nonumber
\\
\leq&3C_2 C_{10}M_2  M_5\varepsilon^2(1+t^{-\frac12})e^{-\alpha_1 t}.\nonumber
\end{align}}
From \eqref{3.60}--\eqref{3.62}, it follows that
\begin{align*}
\|\nabla c\|_{L^\infty(\Omega)}
&\leq (C_3+2 C_2C_{10} M_4 +3C_2 C_{10}M_2  M_5\varepsilon)
(1+t^{-\frac12})\varepsilon e^{-\alpha_1t}
\\
&\leq \frac{M_2}{2}(1+t^{-\frac12})\varepsilon e^{-\alpha_1 t},
\end{align*}
due to the choice of $M_1,M_2$ and $\varepsilon$ satisfying \eqref{33.1}, \eqref{33.2}, and thereby complete the proof.

\begin{lemma}\label{Lemma3.16}
Under the assumptions of Proposition 3.1, for all $\theta\in[q_0,\infty]$ and $t\in(0,T)$,
\begin{align*}
\|(\rho-m)(\cdot,t)-e^{t\Delta}(\rho_0-m_0)\|_{L^\theta(\Omega)}\leq \frac{M_1}{2}\varepsilon (1+t^{-\frac{N}{2}(\frac{1}{p_0}-\frac{1}{\theta})}) e^{-\alpha_1 t}.
\end{align*}
\end{lemma}
\proof According to \eqref{3.53} and Lemma 1.3 of \cite{Winkler7}, we have
{\setlength\abovedisplayskip{4pt}
\setlength\belowdisplayskip{4pt}\begin{align*}
&\|(\rho-m)(\cdot,t)-e^{t\Delta}(\rho_0-m_0)\|_{L^\theta(\Omega)}
\\
\leq&\int_0^t\|e^{(t-s)\Delta}(\nabla\cdot(\rho \mathcal{S}(x,\rho,c)\nabla c)+u\cdot\nabla(\rho-m))(\cdot,s)\|_{L^\theta(\Omega)}ds
\\
\leq&\int_0^t\|e^{(t-s)\Delta}\nabla\cdot(\rho\mathcal{S}(x,\rho,c)\nabla c)(\cdot,s)\|_{L^\theta(\Omega)}ds
+\int_0^t\|e^{(t-s)\Delta}\nabla\cdot((\rho-m-\rho_\infty)u)(\cdot,s)\|_{L^\theta(\Omega)}ds
\\
\leq&C_4C_S\int_0^t(1+(t-s)^{-\frac12-\frac{N}{2}(\frac{1}{q_0}-\frac{1}{\theta})})
e^{-\lambda_1(t-s)}\|\rho(\cdot,s)\|_{L^{q_0}(\Omega)}\|\nabla c(\cdot,s)\|_{L^\infty(\Omega)}ds
\\
&+C_4\int_0^t(1+(t-s)^{-\frac12-\frac{N}{2}(\frac{1}{q_0}-\frac{1}{\theta})})
e^{-\lambda_1(t-s)}\|u(\rho-m-\rho_\infty)(\cdot,s)\|_{L^{q_0}(\Omega)}ds
\\
:=&I_1+I_2.
\end{align*}}
Now we need to estimate $I_1$ and $I_2$. Firstly, from Lemma \ref{Lemma3.11} and Lemma \ref{Lemma3.12}, we obtain
\begin{align}
\|\rho(\cdot,s)\|_{L^{q_0}(\Omega)}&\leq \|(\rho-m-\rho_\infty)(\cdot,s)\|_{L^{q_0}(\Omega)}+\|m(\cdot,s)\|_{L^{q_0}(\Omega)}+\|\rho_\infty\|_{L^{q_0}(\Omega)}\label{3.63}
\\
&\leq M_3\varepsilon
(1+s^{-\frac{N}{2}\left(\frac{1}{p_0}-\frac{1}{q_0}\right)})e^{-\alpha_1s}+M_6 \nonumber
\end{align}
with $M_6=e^{(1+C_1+C_1|\Omega|^{\frac1{p_0}-\frac1{q_0}})\sigma}+\rho_\infty|\Omega|^{\frac{1}{q_0}}$,
which together with Lemma \ref{lemma3.15} and Lemma 1.2 of \cite{Winkler7} implies that\vspace{-0.5em}
{\setlength\abovedisplayskip{4pt}
\setlength\belowdisplayskip{4pt}\begin{align}
I_1&\leq C_4C_SM_6\int_0^t(1+(t-s)^{-\frac12-\frac{N}{2}(\frac{1}{q_0}-\frac{1}{\theta})})
e^{-\lambda_1(t-s)}\|\nabla c\|_{L^\infty(\Omega)}ds
\label{1.52}\\
 &+ M_7\varepsilon
\int_0^t(1+(t-s)^{-\frac12-\frac{N}{2}(\frac{1}{q_0}-\frac{1}{\theta})})(1+s^{-\frac{N}{2}(\frac{1}{p_0}-\frac{1}{q_0})})
e^{-\alpha_1s}
e^{-\lambda_1(t-s)}\|\nabla c\|_{L^\infty(\Omega)}ds
 \nonumber\\
&\leq C_4C_SM_6M_2\varepsilon \int_0^t(1+(t-s)^{-\frac12-\frac{N}{2}(\frac{1}{q_0}-\frac{1}{\theta})})
e^{-\lambda_1(t-s)}(1+s^{-\frac12}) e^{-\alpha_1 s}ds \nonumber\\
& +3 M_7M_2\varepsilon^2
\int_0^t(1+(t-s)^{-\frac12-\frac{N}{2}(\frac{1}{q_0}-\frac{1}{\theta})})(1+s^{-\frac{1}{2}-\frac{N}{2}(\frac{1}{p_0}-\frac{1}{q_0})})
e^{-2\alpha_1s}
e^{-\lambda_1(t-s)}ds\nonumber
\\
&\leq  C_{10} (C_4C_SM_6M_2+3 M_7M_2\varepsilon)(1+t^{-\frac{N}{2}(\frac{1}{p_0}-\frac{1}{\theta})})\varepsilon e^{-\alpha_1 t}\nonumber
\\
&\leq \frac{M_1}{4}\varepsilon(1+t^{-\frac{N}{2}(\frac{1}{p_0}-\frac{1}{\theta})})e^{-\alpha_1 t}\nonumber
\end{align}}
with $M_7=C_4C_SM_3$, where we have used \eqref{33.3} and \eqref{33.4} and  $\frac{1}{p_0}-\frac{1}{q_0}<\frac{1}{N}$.
On the other hand,  from  Lemma \ref{Lemma3.11} and Lemma \ref{Lemma 3.14}, it follows that
{\setlength\abovedisplayskip{4pt}
\setlength\belowdisplayskip{4pt}\begin{align}
I_2&=C_4\int_0^t(1+(t-s)^{-\frac12-\frac{N}{2}(\frac{1}{q_0}-\frac{1}{\theta})})
e^{-\lambda_1(t-s)}\|\rho-m-\rho_\infty\|_{L^\infty(\Omega)}\|u\|_{L^{q_0}(\Omega)}ds\nonumber
\\
&\leq 3C_4M_3M_5\varepsilon^{2}\int_0^t(1+(t-s)^{-\frac12-\frac{N}{2}(\frac{1}{q_0}-\frac{1}{\theta})})
e^{-\lambda_1(t-s)}(1+s^{-\frac12+\frac{N }{ 2q_0}-\frac{N}{2p_0}})e^{-(\alpha_1+\alpha_2)s}
ds\nonumber
\\
&\leq 3C_4M_3M_5C_{10} \varepsilon^{2}(1+t^{\min\{0,\frac{N}{2}
(\frac{1}{\theta}-\frac{1}{p_0})\}})e^{-\min\{\lambda_1,\alpha_1+\alpha_2\}t }\nonumber
\\
&\leq \frac{M_1}{4}\varepsilon (1+t^{-\frac{N}{2}(\frac{1}{p_0}-\frac{1}{\theta})})e^{-\alpha_1 t},\label{1.53}
\end{align}}
where we have used \eqref{33.5} and  $\frac{1}{p_0}-\frac{1}{q_0}<\frac{1}{N}$. Hence combining the above inequalities leads to our conclusion immediately.

{\bf Proof of Theorem 1.3 in the case $\mathcal{S}=0$ on $\partial\Omega$, part 1 (Proposition 3.1).}~
First we claim that $T=T_{max}$. In fact, if  $T<T_{max}$, then by Lemma \ref{lemma3.15} and Lemma \ref{Lemma3.16}, we have
$
\|\nabla c(\cdot,t)\|_{L^{\infty}(\Omega)}\leq \frac{M_2}{2}\varepsilon(1+t^{-\frac12}) e^{-\alpha_1 t}
$
and
{\setlength\abovedisplayskip{4pt}
\setlength\belowdisplayskip{4pt}
\begin{align*}
\|(\rho-m)(\cdot,t)-e^{t\Delta}(\rho_0-m_0)\|_{L^\theta(\Omega)}\leq \frac{M_1}{2}\varepsilon(1+t^{-\frac{N}{2}(\frac{1}{p_0}-\frac{1}{\theta})}) e^{-\alpha_1 t}
\end{align*}}
for all $\theta\in[q_0,\infty]$ and $t\in(0,T)$, which contradicts the definition of $T$ in \eqref{3.52}. Next, we show that $T_{max}=\infty$.
In fact, if $T_{max}<\infty$,
we only need  to show that as  $t \rightarrow T_{max}$,
{\setlength\abovedisplayskip{4pt}
\setlength\belowdisplayskip{4pt}
  $$ \|\rho(\cdot,t)\|_{L^\infty(\Omega)}+\|m(\cdot,t)\|_{L^\infty(\Omega)}
+\|c(\cdot,t)\|_{W^{1,\infty}(\Omega)}+
\|A^{\beta}u(\cdot,t)\|_{L^2(\Omega)}\rightarrow \infty$$
} according to the extensibility criterion in Lemma \ref{lemma2.5}.

Let $t_0:=\min\{1,\frac{T_{max}}3\}$. Then from Lemma \ref{Lemma3.12}, there  exists  $K_1>0$ such that for $t\in(t_0,T_{max})$,
{\setlength\abovedisplayskip{4pt}
\setlength\belowdisplayskip{4pt}\begin{align}\label{3.66a}
\|m(\cdot,t)\|_{L^\infty(\Omega)}\leq K_1e^{-\rho_\infty t}.
\end{align}}
Moreover, from Lemma \ref{Lemma3.11} and the fact that
{\setlength\abovedisplayskip{4pt}
\setlength\belowdisplayskip{4pt}
$$\|\rho(\cdot,t)-\rho_\infty\|_{L^\infty(\Omega)}\leq\|(\rho-m)(\cdot,t)-\rho_\infty\|_{L^\infty(\Omega)}+\|m(\cdot,t)\|_{L^\infty(\Omega)},
$$}
it follows that for  all $t\in(t_0,T_{max})$ and  some constant $K_2>0$,
{\setlength\abovedisplayskip{4pt}
\setlength\belowdisplayskip{4pt}\begin{align}\label{3.67a}
\|\rho(\cdot,t)-\rho_\infty\|_{L^\infty(\Omega)}\leq K_2e^{-\alpha_1 t}.
\end{align}}
Furthermore, Lemma \ref{lemma3.15} implies that there exists $K_3'>0$ such that
{\setlength\abovedisplayskip{4pt}
\setlength\belowdisplayskip{4pt}\begin{align}\label{3.66b}
\|\nabla c(\cdot,t)\|_{L^{\infty}(\Omega)}\leq K_3'e^{-\alpha_2t}\quad \hbox{for all}\,\,t\in(t_0,T_{max}).
\end{align}}
On the other hand, we can conclude that
 $\|c(\cdot,t)\|_{L^\infty(\Omega)}+
\|A^{\beta}u(\cdot,t)\|_{L^2(\Omega)}\leq C~\hbox{for }~t\in(t_0,T_{max})$.
In fact, we first show that there exists  a constant $M_9>0$ such that
{\setlength\abovedisplayskip{4pt}
\setlength\belowdisplayskip{4pt}\begin{align}\label{3.66}
\|A^\beta u(\cdot,t)\|_{L^2(\Omega)}\leq M_9  e^{-\alpha_2 t}
\end{align}}
for $t_0<t<T_{max}$.
By \eqref{3.56}, we have
{\setlength\abovedisplayskip{4pt}
\setlength\belowdisplayskip{4pt}\begin{align*}
\|A^\beta u(\cdot,t)\|_{L^2(\Omega)}\leq\|A^\beta e^{-tA} u_0\|_{L^2(\Omega)}
+\int_{0}^t\|A^\beta e^{-(t-s)A}\mathcal{P}((\rho+m-\rho_\infty)\nabla\Phi)(\cdot,s)\|_{L^2(\Omega)}ds.
\end{align*}}
According to Lemma 2.3 of \cite{Cao1},
$
\|A^\beta e^{-tA} u_0\|_{L^2(\Omega)}\leq  C_5 e^{-\mu t}\|A^\beta u_0\|_{L^2(\Omega)}
$
for all $t\in(0,T_{max})$.
On the other hand, from Lemma 2.3 of \cite{Cao1}  and Lemma \ref{Lemma3.11}, it follows that there exists $\hat{M}>1$ such that
\begin{align*}
&\int_0^t\|A^\beta e^{-(t-s)A}\mathcal{P}((\rho+m-\rho_\infty)\nabla\Phi)(\cdot,s)\|_{L^2(\Omega)}ds
\\
\leq& C_9C_5\|\nabla\Phi\|_{L^\infty(\Omega)}|\Omega|^{\frac{q_0-2}{2q_0}}
\int_0^t e^{-\mu(t-s)}(t-s)^{-\beta} (\|(\rho-m-\rho_\infty)(\cdot,s)\|_{L^{q_0}(\Omega)}+2\|m(\cdot,s)\|_{L^{q_0}(\Omega)})ds
\\
\leq& C_9C_5\|\nabla\Phi\|_{L^\infty(\Omega)}|\Omega|^{\frac{q_0-2}{2q_0}} \hat{M}
\int_0^t e^{-\mu(t-s)}(t-s)^{-\beta} (1+s^{-\frac{N}{2}(\frac{1}{p_0}-\frac{1}{q_0})})e^{-\alpha_1s}
ds\\
\leq& C_5 C_9 C_{10}\|\nabla\Phi\|_{L^\infty(\Omega)}|\Omega|^{\frac{q_0-2}{2q_0}} \hat{M}
e^{-\alpha_2 t}(1+t^{\min\{0,1-\beta-\frac{N}{2}(\frac{1}{p_0}-\frac{1}{q_0})\}})\\
\leq& C_5 C_9 C_{10}\|\nabla\Phi\|_{L^\infty(\Omega)}|\Omega|^{\frac{q_0-2}{2q_0}} \hat{M}
e^{-\alpha_2 t}(1+t_0^{\min\{0,1-\beta-\frac{N}{2}(\frac{1}{p_0}-\frac{1}{q_0})\}})
\end{align*}
for $t_0<t<T_{max}$. Hence combining the above inequalities, we arrive at \eqref{3.66}.

Since $D(A^\beta)\hookrightarrow L^\infty(\Omega)$ with $\beta\in(\frac{N}{4},1)$, we have
{\setlength\abovedisplayskip{4pt}
\setlength\belowdisplayskip{4pt}
\begin{align}\label{3.67}
\|u(\cdot,t)\|_{L^\infty(\Omega)}\leq K_4  e^{-\alpha_2 t} \quad\hbox{for some}\,K_4>0 \,\hbox{and }\,  t\in(0,T_{max}).
\end{align}}
Now we turn to show that
 there exists $K_3''>0$ such that
 {\setlength\abovedisplayskip{4pt}
\setlength\belowdisplayskip{4pt}
\begin{align}\label{3.71}
\|c(\cdot,t)\|_{L^{\infty}(\Omega)}\leq K_3'' e^{-\alpha_2t}\quad\hbox{for all}\,\, t\in(0,T_{max}).
\end{align}}
Indeed, from \eqref{3.55}, it follows that
\begin{align}
\|c\|_{L^\infty(\Omega)}
&\leq\|e^{t(\Delta-1)}c_0\|_{L^\infty(\Omega)}+\int_0^t\|
e^{(t-s)(\Delta-1)}(m-u\cdot\nabla c)\|_{L^\infty(\Omega)}ds\nonumber
\\
&\leq  e^{-t}\|c_0\|_{L^\infty(\Omega)}+\int_0^t\|
e^{(t-s)(\Delta-1)}m(\cdot,s)\|_{L^\infty(\Omega)}ds\label{3.72}
\\
&\quad+\int_0^t\|
e^{(t-s)(\Delta-1)}u\cdot\nabla c(\cdot,s)\|_{L^\infty(\Omega)}ds.\nonumber
\end{align}
 An application of \eqref{3.57} with $k=\infty$ yields
{\setlength\abovedisplayskip{4pt}
\setlength\belowdisplayskip{4pt}\begin{align}
\int_0^t\|e^{(t-s)(\Delta-1)}m(\cdot,s)\|_{L^\infty(\Omega)}ds
&\leq \int_0^t e^{-(t-s)}\|m(\cdot,s)\|_{L^{\infty}(\Omega)}ds\label{3.73}
\\
&\leq \|m_0\|_{L^\infty(\Omega)}M_4\int_0^t e^{-(t-s)} e^{-\rho_\infty s}
ds\nonumber
\\
& \leq M_4C_{10} e^{-\alpha_2t}.\nonumber
\end{align}}
On the other hand,
 from  \eqref{3.67} and \eqref{3.66b}, we can see  that
 {\setlength\abovedisplayskip{4pt}
\setlength\belowdisplayskip{4pt}
\begin{align}\label{3.74}
\int_0^t\|
e^{(t-s)(\Delta-1)} u\cdot\nabla c\|_{L^\infty(\Omega)}ds
&\leq \int_0^t e^{-(t-s)}\|u\|_{L^{\infty}(\Omega)}\|\nabla c\|_{L^\infty(\Omega)}ds
\\
&\leq K_3' K_4  \int_0^t e^{-2\alpha_2 s} e^{-(t-s)}ds
\nonumber
\\ &\leq K_3' K_4 C_{10} e^{-\alpha_2 t}.
 \nonumber
\end{align}}
Hence, inserting  \eqref{3.73}, \eqref{3.74} into  \eqref{3.72}, we arrive at the conclusion \eqref{3.71}. Therefore
we have $T_{max}=\infty$, and the decay estimates in \eqref{3.48}--\eqref{3.51} follow from \eqref{3.66a}--\eqref{3.71}, respectively. 
\vspace{-1em}

\subsubsection{The case $\int_{\Omega}\rho_0<\int_{\Omega}m_0$}\vspace{-0.5em}

In this subsection we consider the case $\int_{\Omega}\rho_0<\int_{\Omega}m_0$, i.e.,
$m_\infty>0$, $\rho_\infty=0$.\vspace{-0.5em}
\begin{proposition}
Suppose that \eqref{1.4} hold with $\alpha=0$ and $\int_{\Omega}\rho_0<\int_{\Omega}m_0$. Let $N=3$, $p_0\in(\frac {2N}3, N)$, $q_0\in(N,\frac{Np_0}{2(N-p_0)})$. 
Then there exists
$\varepsilon>0$ such that for any initial data $(\rho_0,m_0,c_0,u_0)$ fulfilling
\eqref{1.7} as well as
{\setlength\abovedisplayskip{4pt}
\setlength\belowdisplayskip{4pt}
$$\|\rho_0\|_{L^{p_0}(\Omega)}\leq\varepsilon,\quad  \|m_0-m_\infty\|_{L^{q_0}(\Omega)}\leq\varepsilon,
\quad\|\nabla c_0\|_{L^{N}(\Omega)}\leq\varepsilon, \quad\|u_0\|_{L^{N}(\Omega)}\leq\varepsilon,
$$}
\eqref{1.1} admits a global classical solution $(\rho,m,c,u,P)$.
Furthermore,  for any $\alpha_1\!\in\!(0,\min\{\lambda_1,m_\infty\})$, $\alpha_2\!\in\!(0,\min\{\alpha_1,\lambda_1',1\})$, there exist constants $K_i>0$, $i=1,2,3,4$, such that
\begin{align}
&\|m(\cdot,t)-m_\infty\|_{L^\infty(\Omega)}\leq K_1e^{-\alpha_1 t},\label{3.75}
\\
&\|\rho(\cdot,t)\|_{L^\infty(\Omega)}\leq K_2e^{-\alpha_1 t},\label{3.76}
\\
&\|c(\cdot,t)-m_\infty\|_{W^{1,\infty}(\Omega)}\leq K_3e^{-\alpha_2t},\label{3.77}
\\
&\|u(\cdot,t)\|_{L^\infty(\Omega)}\leq K_4 e^{-\alpha_2 t}.\label{3.78}
\end{align}
\end{proposition}

The proof of Proposition 3.2 proceeds in a parallel fashion to that of Proposition 3.1.  
However, due to differences in the properties of $\rho$ and $m$, 
there are significant differences in the details of their proofs. Thus, for the convenience of the reader, we will give the full proof of Proposition 3.2.

The following can be verified easily:\vspace{-1em}
\begin{lemma}\label{Lemma3.17a}
Under the assumptions of Proposition 3.2, 
 it is possible to choose $M_1>0,M_2>0$ and $\varepsilon>0$ such that
\begin{align}
&C_3\leq \displaystyle \frac {M_2}6,\quad
C_2 C_{10}(1+ C_1+C_1|\Omega|^{\frac1{p_0}-\frac1{q_0}}+M_1)\leq \displaystyle \frac {M_2}6,\label{34.1}
\\
&18C_2C_6 C_{10} (1+2C_9 C_{10} (1+C_1+ C_1|\Omega|^{\frac1{p_0}-\frac1{q_0}}+2M_1)
\|\nabla\Phi\|_{L^\infty(\Omega)}) \varepsilon \leq 1, \label{34.2}
\\
&2C_1+(\min\{1,|\Omega|\})^{-\frac 1{p_0}}\leq \displaystyle\frac {M_1}8,\quad
24C_4C_SC_{10}M_2 \varepsilon<1,\label{34.3}
\\
&24C_4C_{10} C_6
( 1+2C_9 C_{10}
(1+C_1+ C_1|\Omega|^{\frac1{p_0}-\frac1{q_0}}+2M_1)
  \|\nabla\Phi\|_{L^\infty(\Omega)})
\varepsilon<1,\label{34.4}
\\
&24C_4C_{10}(1+ C_1+C_1|\Omega|^{\frac1{p_0}-\frac1{q_0}}+M_1)\varepsilon<1,\label{34.5}
\\
&12C_4C_SC_{10} M_1M_2 \varepsilon< 1,\label{34.6}
\\
&C_{10}C_6 C_4(1\!+\! C_1\!+\!C_1|\Omega|^{\frac1{p_0}\!-\!\frac1{q_0}})
(1\!+\!2C_9 C_{10} (1\!+\!C_1\!+\! C_1|\Omega|^{\frac1{p_0}\!-\!\frac1{q_0}}\!+\!2M_1)
\|\nabla\Phi\|_{L^\infty(\Omega)})\varepsilon \!< \!\displaystyle\frac1 {24}.\label{34.7}
\end{align}
\end{lemma}

Similar to the proof of  Proposition 3.1, we define{\setlength\abovedisplayskip{4pt}
\setlength\belowdisplayskip{4pt}
\begin{align}
T\!\triangleq\!\sup\!\left\{\!\widetilde{T}\!\in\!(0,T_{max})\!\left|
\begin{array}{ll}
\!\|(m\!-\!\rho)(\cdot,t)\!-\!e^{t\Delta}(m_0\!-\!\rho_0)\|_{L^{\theta}(\Omega)}\!\leq\! \varepsilon(1+t^{-\frac{N}{2}(\frac{1}{p_0}-\frac{1}{\theta})})e^{-\alpha_1t},
\\
\!\|\rho(\cdot,t)\|_{L^\theta(\Omega)}\leq M_1\varepsilon(1+t^{-\frac{N}{2}(\frac{1}{p_0}-\frac{1}{\theta})})e^{-\alpha_1t},\,\forall\theta\in[q_0,\infty],
\\
\!\|\nabla c(\cdot,t)\|_{L^\infty(\Omega)}\leq M_2\varepsilon(1+t^{-\frac12})e^{-\alpha_1 t} \quad \hbox{for all }\, t\in[0,\widetilde{T}).
\end{array}\right.
\!\right\}\label{3.79}
\end{align}}
By Lemma 2.1 and (1.7), $T>0$ is well-defined. As in the previous subsection,  we first show  $T=T_{max}$, and then $T_{max}=\infty$. To this end,
we will show that all of the estimates
mentioned in \eqref{3.79} are valid with even smaller coefficients on the right hand side than appearing in \eqref{3.79}. The derivation of these estimates will mainly rely on $L^p-L^q$ estimates for the Neumann heat semigroup and the corresponding semigroup for Stokes operator, and the fact that the classical solutions of \eqref{1.1} on $(0,T)$ can be represented as{\setlength\abovedisplayskip{4pt}
\setlength\belowdisplayskip{4pt}
\begin{align}
&(m-\rho)(\cdot,t)\!=\!e^{t\Delta}(m_0-\rho_0)\!+\!\int_0^te^{(t-s)\Delta}(\nabla\cdot(\rho \mathcal{S}(x,\rho,c)\nabla c)-u\cdot\nabla(m-\rho))(\cdot,s)ds,\label{3.80}
\\
&\rho(\cdot,t)=e^{t\Delta}\rho_0-\int_0^te^{(t-s)\Delta}(\nabla\cdot(\rho \mathcal{S}(x,\rho,c)\nabla c)-u\cdot\nabla\rho+\rho m)(\cdot,s)ds,\label{3.81}
\\
&c(\cdot,t)=e^{t(\Delta-1)}c_0+\int_0^te^{(t-s)(\Delta-1)}(m-u\cdot\nabla c)(\cdot,s)ds,\label{3.82}
\\
&u(\cdot,t)=e^{-tA}u_0+\int_0^te^{-(t-s)A}\mathcal{P}((\rho+m)\nabla\Phi)(\cdot,s)ds.\label{3.83}
\end{align}}
\begin{lemma}\label{lemma3.18}
Under the assumptions of Proposition 3.2, we have
\begin{align*}
\|(m-\rho)(\cdot,t)-m_\infty\|_{L^\theta(\Omega)}\leq M_3\varepsilon(1+t^{-\frac{N}{2}(\frac{1}{p_0}-\frac{1}{\theta})})e^{-\alpha_1t}
\end{align*}
for all $t\in(0,T)$ and $\theta\in[q_0,\infty]$.
\end{lemma}

\proof  Since $e^{t\Delta}(\overline{m}_0-\overline{\rho}_0)=m_\infty$ and $\int_\Omega (m_0-\rho_0-m_\infty)=0$, from the Definition of $T$ and Lemma 1.3 of \cite{Winkler7}, we get{\setlength\abovedisplayskip{4pt}
\setlength\belowdisplayskip{4pt}
\begin{align*}
&\|(m-\rho)(\cdot,t)-m_\infty\|_{L^\theta(\Omega)}
\\
\leq&\|(m-\rho)(\cdot,t)-e^{t\Delta}(m_0-\rho_0)\|_{L^\theta(\Omega)}
+\|e^{t\Delta}(m_0-\rho_0)-e^{t\Delta}m_{\infty}\|_{L^\theta(\Omega)}
\\
\leq& \varepsilon(1+t^{-\frac{N}{2}(\frac{1}{p_0}-\frac{1}{\theta})})e^{-\alpha_1t}
+C_1(1+t^{-\frac{N}{2}(\frac{1}{p_0}-\frac{1}{\theta})})(\|\rho_0\|_{L^{p_0}(\Omega)}
+\|m_0-m_\infty\|_{L^{p_0}(\Omega)})e^{-\lambda_1t}
\\
\leq&(1+  C_1+C_1|\Omega|^{\frac1{p_0}-\frac1{q_0}})\varepsilon(1+t^{-\frac{N}{2}(\frac{1}{p_0}-\frac{1}{\theta})})e^{-\alpha_1t}
\end{align*}}
for all $t\in(0,T)$ and $\theta\in[q_0,\infty]$. This lemma is proved for $M_3=1+  C_1+C_1|\Omega|^{\frac1{p_0}-\frac1{q_0}}$.

\begin{lemma}\label{lemma3.19}
Under the assumptions of Proposition 3.2, we have
{\setlength\abovedisplayskip{3pt}
\setlength\belowdisplayskip{3pt}
\begin{align*}
\|m(\cdot,t)-m_\infty\|_{L^{\theta}(\Omega)}\leq M_4\varepsilon(1+t^{-\frac N2(\frac{1}{p_0}-\frac{1}{\theta})}) e^{-\alpha_1 t}
\quad\hbox{for all}\,\,t\in(0,T), \theta\in[q_0,\infty].
\end{align*}}
\end{lemma}
\proof From Lemma \ref{lemma3.18} and the definition of $T$, it follows that
{\setlength\abovedisplayskip{4pt}
\setlength\belowdisplayskip{4pt}
\begin{align*}
\|m(\cdot,t)-m_\infty\|_{L^{\theta}(\Omega)}
\leq &\|(m-\rho-m_\infty)(\cdot,t)\|_{L^{\theta}(\Omega)}+\|\rho(\cdot,t)\|_{L^{\theta}(\Omega)}
\\
\leq& (M_3+M_1)\varepsilon(1+t^{-\frac N2(\frac{1}{p_0}-\frac{1}{\theta})})e^{-\alpha_1 t}.
\end{align*}}
The Lemma is proved for $M_4=M_3+M_1$.

\begin{lemma}\label{lemma3.20}
Under the assumptions of Proposition 3.2, there exists $M_5>0$ such that
\begin{align*}
\|u(\cdot,t)\|_{L^{q_0}(\Omega)}\leq M_5\varepsilon(1+t^{-\frac12+\frac{N}{2q_0}}) e^{-\alpha_2 t}
\quad ~\hbox{for all }~t\in(0,T).\end{align*}
\end{lemma}

\proof For any given $\alpha_2<\lambda_1'$, we can fix $ \mu\in (\alpha_2, \lambda_1')$.
 By  \eqref{3.83}, Lemma 2.3 of \cite{Cao1}, and noticing that  $\mathcal{P}(\nabla \Phi)=0$, we obtain that
{\setlength\abovedisplayskip{4pt}
\setlength\belowdisplayskip{4pt}\begin{align}
&\|u(\cdot,t)\|_{L^{q_0}(\Omega)}\nonumber
\\
\leq& C_6t^{-\frac N2(\frac1N-\frac{1}{q_0})}e^{-\mu t}\|u_0\|_{L^N(\Omega)}
+\int_0^t\|e^{-(t-s)A}\mathcal{P}((\rho+m)\nabla\Phi)(\cdot,s)\|_{L^{q_0}(\Omega)}ds\label{3.84}
\\
\leq& C_6t^{-\frac N2(\frac 1N-\frac{1}{q_0})}e^{-\mu t}\|u_0\|_{L^N(\Omega)}
+C_6C_9\int_0^te^{-\mu(t-s)}\|(\rho+m-m_\infty)(\cdot,s)\|_{L^{q_0}(\Omega)}\|\nabla\Phi\|_{L^\infty(\Omega)}ds\nonumber
\\
\leq& C_6 t^{-\frac 12+\frac{N}{2q_0}}e^{-\mu t}\|u_0\|_{L^N(\Omega)}
+C_6C_9\|\nabla\Phi\|_{L^\infty(\Omega)}\int_0^te^{-\mu(t-s)}\|(\rho+m-m_\infty)(\cdot,s)\|_{L^{q_0}(\Omega)}ds.\nonumber
\end{align}}
By Lemma \ref{lemma3.19} and the definition of $T$, we get
{\setlength\abovedisplayskip{4pt}
\setlength\belowdisplayskip{4pt}\begin{align}
\|(\rho+m-m_\infty)(\cdot,s)\|_{L^{q_0}(\Omega)}
=&\|(m-m_\infty)(\cdot,s)\|_{L^{q_0}(\Omega)}+\|\rho(\cdot,s)\|_{L^{q_0}(\Omega)}\label{3.85}
\\
\leq &(M_4+M_1)\varepsilon(1+s^{-\frac{N}{2}(\frac{1}{p_0}-\frac{1}{q_0})})e^{-\alpha_1s}.\nonumber
\end{align}}
Inserting \eqref{3.85} into \eqref{3.84},
and noting  $\frac{N}{2}(\frac{1}{p_0}-\frac{1}{q_0})<1$, we have 
{\setlength\abovedisplayskip{4pt}
\setlength\belowdisplayskip{4pt}\begin{align*}
&\|u(\cdot,t)\|_{L^{q_0}(\Omega)}
\\
\leq& C_6t^{-\frac 12+\frac{N}{2q_0}}e^{-\mu t}\|u_0\|_{L^N(\Omega)}
+C_6C_9(M_4+M_1)\|\nabla\Phi\|_{L^\infty(\Omega)}\varepsilon\int_0^t(1+s^{-\frac N2(\frac{1}{p_0}-\frac{1}{q_0})})e^{-\alpha_1s}e^{-\mu(t-s)}ds
\\
\leq& C_6 t^{-\frac12+\frac{N}{2q_0}}e^{-\mu t}\|u_0\|_{L^N(\Omega)}
+C_6C_9 C_{10}(M_4+M_1)\|\nabla\Phi\|_{L^\infty(\Omega)}\varepsilon(1+t^{\min\{0,1-\frac N2(\frac{1}{p_0}-\frac{1}{q_0})\}})e^{-\alpha_2t}
\\
\leq& C_6 t^{-\frac12+\frac{N}{2q_0}}\varepsilon e^{-\mu t}
+2C_6C_9 C_{10}(M_4+M_1) \|\nabla\Phi\|_{L^\infty(\Omega)}\varepsilon e^{-\alpha_2 t}
\\
= & M_5\varepsilon(1+t^{-\frac12+\frac{N}{2q_0}})e^{-\alpha_2 t}
\end{align*}}
with  $M_5=C_6+2C_6C_9 C_{10}(M_4+M_1)\|\nabla\Phi\|_{L^\infty(\Omega)}$. \vspace{-0.5em}

\begin{lemma}\label{lemma3.21}
Under the assumptions of Proposition  3.2, we have{\setlength\abovedisplayskip{4pt}
\setlength\belowdisplayskip{4pt}
\begin{align*}
\|\nabla c(\cdot,t)\|_{L^{\infty}(\Omega)}\leq \frac{M_2}{2}\varepsilon(1+t^{-\frac12}) e^{-\alpha_1 t}\quad \hbox{for all}\,\, t\!\in\!(0,T).
\end{align*}}
\end{lemma}
\vspace{-1em}
\proof From \eqref{3.82} and Lemma 1.3 of \cite{Winkler7}, we have{\setlength\abovedisplayskip{4pt}
\setlength\belowdisplayskip{4pt}
\begin{align}
\|\nabla c(\cdot,t)\|_{L^\infty(\Omega)}
&\leq\|e^{t(\Delta-1)}\nabla c_0\|_{L^\infty(\Omega)}+\int_0^t\|\nabla
e^{(t-s)(\Delta-1)}(m-u\cdot\nabla c)(\cdot,s)\|_{L^\infty(\Omega)}ds\nonumber
\\
&\leq C_3(1+t^{-\frac12})e^{-(\lambda_1+1)t}\|\nabla c_0\|_{L^N(\Omega)}+\int_0^t\|\nabla
e^{(t-s)(\Delta-1)}(m-m_{\infty})(\cdot,s)\|_{L^\infty(\Omega)}ds\nonumber
\\
&\quad+\int_0^t\|\nabla
e^{(t-s)(\Delta-1)}u\cdot\nabla c(\cdot,s)\|_{L^\infty(\Omega)}ds.\label{3.86}
\end{align}}
In the second inequality, we have used $ \nabla e^{(t-s)(\Delta-1)}m_{\infty}=0$.

 From Lemma 1.3 of \cite{Winkler7}, Lemma \ref{lemma3.19}, Lemma 1.2 of \cite{Winkler7},  it follows that
{\setlength\abovedisplayskip{4pt}
\setlength\belowdisplayskip{4pt}
\begin{align}
&\int_0^t\|\nabla
e^{(t-s)(\Delta-1)}(m-m_{\infty})(\cdot,s)\|_{L^\infty(\Omega)}ds\nonumber
\\
\leq&C_2\int_0^t(1+(t-s)^{-\frac12-\frac{N}{2q_0}})e^{-(\lambda_1+1)(t-s)}\|(m-m_{\infty}) (\cdot,s)\|_{L^{q_0}(\Omega)}ds\label{3.87}
\\
\leq&C_2 M_4\varepsilon
\int_0^t(1+(t-s)^{-\frac12-\frac{N}{2q_0}})e^{-(\lambda_1+1)(t-s)}
(1+s^{-\frac N2(\frac{1}{p_0}-\frac{1}{q_0})}) e^{-\alpha_1 s}ds\nonumber
\\
\leq&C_2 C_{10}M_4\varepsilon (1+t^{\min\{0,\frac12-\frac{N}{2p_0}\}})e^{-\min\{\alpha_1,\lambda_1+1\}t}\nonumber
\\
\leq&C_2 C_{10}M_4\varepsilon
(1+t^{-\frac12})e^{-\alpha_1t}.\nonumber
\end{align}}
On the other hand,  by Lemma 1.3 of \cite{Winkler7},  Lemma \ref{lemma3.20} and the definition of $T$, we obtain
{\setlength\abovedisplayskip{4pt}
\setlength\belowdisplayskip{4pt}
\begin{align}
&\int_0^t\|\nabla
e^{(t-s)(\Delta-1)}u\cdot\nabla c(\cdot,s)\|_{L^\infty(\Omega)}ds\nonumber
\\
\leq&C_2\int_0^t(1+(t-s)^{-\frac12-\frac{N}{2q_0}})e^{-(\lambda_1+1)(t-s)}\|u\cdot\nabla c(\cdot,s)\|_{L^{q_0}(\Omega)}ds\label{3.88}
\\
\leq&C_2\int_0^t(1+(t-s)^{-\frac12-\frac{N}{2q_0}})e^{-(\lambda_1+1)(t-s)}\|u(\cdot,s)\|_{L^{q_0}(\Omega)}\|\nabla c(\cdot,s)
\|_{L^\infty(\Omega)}ds\nonumber
\\
\leq& C_2M_5M_2\varepsilon^2
\int_0^t(1+(t-s)^{-\frac12-\frac{N}{2q_0}})e^{-(\lambda_1+1)(t-s)}
(1+s^{-\frac12+\frac{N}{2q_0}}) (1+s^{-\frac12}) e^{-(\alpha_1+\alpha_2) s} \nonumber
\\
\leq&3C_2M_5M_2\varepsilon^{2}\int_0^te^{-(\lambda_1+1)(t-s)}
 e^{-(\alpha_1+\alpha_2) s}(1+(t-s)^{-\frac12-\frac{N}{2q_0}})(1+s^{-1+\frac{N}{2q_0}})ds\nonumber
\\
\leq&3C_2M_5M_2C_{10}\varepsilon^{2}(1+t^{-\frac12})e^{-\min\{\lambda_1+1,\alpha_1+\alpha_2\}t}\nonumber
\\
\leq&3C_2M_5M_2C_{10}\varepsilon^{2}(1+t^{-\frac12})e^{-\alpha_1t}. \nonumber
\end{align}}
Hence combining above inequalities with \eqref{34.1}, \eqref{34.2}, we arrive at the conclusion.\vspace{-0.5em}

\begin{lemma}\label{lemma3.22}
Under the assumptions of Proposition 3.2, we have
\begin{align*}
\|\rho(\cdot,t)\|_{L^\theta(\Omega)}\leq \frac{M_1}{2}\varepsilon(1+t^{-\frac{N}{2}(\frac{1}{p_0}-\frac{1}{\theta})})e^{-\alpha_1t}
\quad\hbox{for all}\,\,t\in(0,T),\,\theta\in[q_0,\infty].
\end{align*}
\end{lemma}
\proof By the variation-of-constants formula, we have 
\begin{align*}
\rho(\cdot,t)=&e^{t(\Delta-m_\infty)}\rho_0-\int_0^te^{(t-s)(\Delta-m_\infty)}(\nabla\cdot(\rho\mathcal{S}(\cdot,\rho,c)\nabla c)-u\cdot\nabla\rho)(\cdot,s)ds\\
&+
\int_0^te^{(t-s)(\Delta-m_\infty)} \rho (m_\infty-m)(\cdot,s)ds
.
\end{align*}
By Lemma 1.3 of \cite{Winkler7}, 
the result in Section 2 of \cite{Horstmann 2} and noticing $\alpha_1<\min\{\lambda_1,m_\infty\}$, we obtain
{\setlength\abovedisplayskip{4pt}
\setlength\belowdisplayskip{4pt}
\begin{align*}
&\|\rho(\cdot,t)\|_{L^\theta(\Omega)}
\\
\leq&
e^{-m_\infty t}(\|e^{t\Delta}(\rho_0-\overline{\rho}_0)\|_{L^\theta(\Omega)}+ \|\overline{\rho}_0 \|_{L^\theta(\Omega)}) +\int_0^t\|e^{(t-s)(\Delta-m_\infty)}\nabla\cdot(\rho \mathcal{S}(\cdot,\rho,c)\nabla c)(\cdot,s)\|_{L^\theta(\Omega)}ds
\\
&+\int_0^t\|e^{(t-s)(\Delta-m_\infty)}(u\cdot\nabla\rho)(\cdot,s)\|_{L^\theta(\Omega)}ds
+\int_0^t\|e^{(t-s)(\Delta-m_\infty)}\rho (m_\infty-m)(\cdot,s)\|_{L^\theta(\Omega)}ds
\\
\leq&
C_1(1+t^{-\frac N2(\frac1{p_0}-\frac1{\theta})})e^{-(\lambda_1+m_\infty)t}\|\rho_0-\overline{\rho}_0\|_{L^{p_0}(\Omega)}
+(\min\{1,|\Omega|\})^{-\frac 1{p_0}}e^{-m_\infty t}\varepsilon \\
& +C_4C_S\int_0^t(1+(t-s)^{-\frac12-\frac{N}{2}(\frac{1}{q_0}-\frac{1}{\theta})})e^{-(\lambda_1+m_\infty)(t-s)}\|\rho\|_{L^{q_0}(\Omega)}\|\nabla c\|_{L^\infty(\Omega)}ds\\
&+\int_0^t\|e^{(t-s)(\Delta-m_\infty) }\nabla\cdot(\rho u)(\cdot,s)\|_{L^\theta(\Omega)}ds
+\int_0^t\|e^{(t-s)(\Delta-m_\infty)}\rho (m_\infty-m)(\cdot,s)\|_{L^\theta(\Omega)}ds\\
\leq&
(2C_1+(\min\{1,|\Omega|\})^{-\frac 1{p_0}}) (1+t^{-\frac N2(\frac1{p_0}-\frac1{\theta})}) \varepsilon e^{-\alpha_1t}
\\
&+C_4C_S\int_0^t(1+(t-s)^{-\frac12-\frac{N}{2}(\frac{1}{q_0}-\frac{1}{\theta})})
e^{-(\lambda_1+m_\infty)(t-s)}\|\rho\|_{L^{q_0}(\Omega)}\|\nabla c\|_{L^\infty(\Omega)}ds
\\
&+C_4\int_0^t(1+(t-s)^{-\frac12-\frac{N}{2}(\frac{1}{q_0}-\frac{1}{\theta})})
e^{-(\lambda_1+m_\infty)(t-s)}\|\rho\|_{L^{\infty}(\Omega)}\|u\|_{L^{q_0}(\Omega)}ds
\\
& +C_1\int_0^t(1+(t-s)^{-\frac{N}{2}(\frac{1}{q_0}-\frac{1}{\theta})})
e^{-m_\infty(t-s)}\|\rho\|_{L^{q_0}(\Omega)}\|m-m_\infty\|_{L^{\infty }(\Omega)}ds
\\
=&(2C_1+(\min\{1,|\Omega|\})^{-\frac 1{p_0}}) (1+t^{-\frac N2(\frac1{p_0}-\frac1{\theta})}) \varepsilon e^{-\alpha_1t}
+I_1+I_2+I_3.
\end{align*}}
By the definition of $T$,  Lemma \ref{lemma3.21}, Lemma 1.2 of \cite{Winkler7} and \eqref{34.3}, we get
\begin{align}
I_1&\leq 3 C_4C_S M_1M_2\varepsilon^2 \int_0^t(1+(t-s)^{-\frac12-\frac{N}{2}(\frac{1}{q_0}-\frac{1}{\theta})})
e^{-\lambda_1(t-s)}e^{-2\alpha_1 s}(1+s^ {-\frac12-\frac{N}{2}(\frac{1}{p_0}-\frac{1}{q_0})})
ds\nonumber
\\
&\leq3C_4C_SC_{10} M_1M_2 \varepsilon^2
(1+t^{\min\{0,-\frac{N}{2}(\frac{1}{p_0}-\frac{1}{\theta})\}})
e^{-\min\{\lambda_1,2\alpha_1\}t}\nonumber
\\
&\leq\frac{M_1}{8}\varepsilon(1+t^{-\frac N2(\frac1{p_0}-\frac1{\theta})})e^{-\alpha_1t}. \nonumber
\end{align}
Similarly, by \eqref{34.5} and \eqref{34.6},  we can also get
{\setlength\abovedisplayskip{4pt}
\setlength\belowdisplayskip{4pt}
\begin{align}
I_2&\leq
3 C_4 M_1M_5\varepsilon^2 \int_0^t(1+(t-s)^{-\frac12-\frac{N}{2}(\frac{1}{q_0}-\frac{1}{\theta})})
e^{-\lambda_1(t-s)}e^{-2\alpha_1 s}(1+s^ {-\frac12-\frac{N}{2}(\frac{1}{p_0}-\frac{1}{q_0})})
ds\nonumber
\\
&\leq3C_4C_{10}M_5M_1\varepsilon^2(1+t^{\min\{0,-\frac{N}{2}(\frac{1}{p_0}-\frac{1}{\theta})\}})
e^{-\min\{\lambda_1,2\alpha_1\}t}\nonumber
\\
&\leq\frac{M_1}{8}\varepsilon(1+t^{-\frac N2(\frac1{p_0}-\frac1{\theta})})e^{-\alpha_1t},\nonumber
\end{align}
\begin{align}
I_3&\leq
3 C_4 M_1M_4\varepsilon^2 \int_0^t(1+(t-s)^{-\frac12-\frac{N}{2}(\frac{1}{q_0}-\frac{1}{\theta})})
e^{-m_\infty(t-s)}
  e^{-2\alpha_1 s}(1+s^{-\frac N{p_0}+\frac{ N}{2 q_0}})
ds\nonumber
\\
&\leq3C_4C_{10}M_1M_4
\varepsilon^2(1+t^{\min\{0,-\frac{N}{2}(\frac{1}{p_0}-\frac{1}{\theta})\}})
e^{-\min\{m_\infty,2\alpha_1\}t}\nonumber
\\
&\leq\frac{M_1}8\varepsilon(1+t^{-\frac N2(\frac1{p_0}-\frac1{\theta})})e^{-\alpha_1t},\nonumber
\end{align}}
respectively, where the fact that $q_0\in (N,\frac{Np_0}{2(N-p_0)})$ warrants $-\frac N{p_0}+\frac{ N}{2 q_0}>-1$ is used.
Hence the combination of the above inequalities  yields
$\|\rho(\cdot,t)\|_{L^\theta(\Omega)}\leq\frac{M_1}{2}\varepsilon(1+t^{-\frac N2(\frac1{p_0}-\frac1{\theta})})e^{-\alpha_1t}.$

\begin{lemma}\label{lemma3.23}
Under the assumptions of Proposition 3.2, we have{\setlength\abovedisplayskip{4pt}
\setlength\belowdisplayskip{4pt}
\begin{align*}
\|(m-\rho)(\cdot,t)-e^{t\Delta}(m_0-\rho_0)\|_{L^\theta(\Omega)}\leq \frac{\varepsilon }{2}(1+t^{-\frac{N}{2}(\frac{1}{p_0}-\frac{1}{\theta})}) e^{-\alpha_1 t}
\,\,\hbox{for }\,\theta\in[q_0,\infty], t\in(0,T).\end{align*}}
\end{lemma}

\proof From \eqref{3.80} and Lemma 1.3 of \cite{Winkler7}, it follows that
{\setlength\abovedisplayskip{4pt}
\setlength\belowdisplayskip{4pt}
\begin{align*}
&\|(m-\rho)(\cdot,t)-e^{t\Delta}(m_0-\rho_0)\|_{L^\theta(\Omega)}
\\
\leq&\int_0^t\|e^{(t-s)\Delta}(\nabla\cdot(\rho \mathcal{S}(\cdot,\rho,c)\nabla c)-u\cdot\nabla(m-\rho))(\cdot,s)\|_{L^\theta(\Omega)}ds
\\
\leq&\int_0^t\|e^{(t-s)\Delta}\nabla\cdot(\rho\mathcal{S}(\cdot,\rho,c)\nabla c)(\cdot,s)\|_{L^\theta(\Omega)}ds
+\int_0^t\|e^{(t-s)\Delta}\nabla\cdot((m-\rho-m_\infty)u)(\cdot,s)\|_{L^\theta(\Omega)}ds
\\
\leq&C_4C_S \int_0^t(1+(t-s)^{-\frac12-\frac{N}{2}(\frac{1}{q_0}-\frac{1}{\theta})})
e^{-\lambda_1(t-s)}\|\rho(\cdot,s)\|_{L^{q_0}(\Omega)}\|\nabla c(\cdot,s)\|_{L^\infty(\Omega)}ds
\\
&+C_4\int_0^t(1+(t-s)^{-\frac12-\frac{N}{2}(\frac{1}{q_0}-\frac{1}{\theta})})
e^{-\lambda_1(t-s)}\|u(m-\rho-m_\infty)(\cdot,s)\|_{L^{q_0}(\Omega)}ds
\\
=&I_1+I_2.
\end{align*}}
\vspace{-0.5em}
From the definition of $T$ and \eqref{34.7}, we have
\begin{align*}
I_1\leq &C_4C_S M_1 M_2 \varepsilon^2 \int_0^t(1+(t-s)^{-\frac12-\frac{N}{2}(\frac{1}{q_0}-\frac{1}{\theta})})
e^{-\lambda_1(t-s)}
(1+s^{-\frac12-\frac{N}{2}(\frac{1}{p_0}-\frac{1}{q_0})})
e^{-2\alpha_1s}ds
\nonumber\\
& \leq3C_4C_SC_{10} M_1M_2 \varepsilon^2
(1+t^{\min\{0,-\frac{N}{2}(\frac{1}{p_0}-\frac{1}{\theta})\}})
e^{-\min\{\lambda_1,2\alpha_1\}t}\nonumber
\\
&\leq
\frac{\varepsilon }{4}(1+t^{-\frac N2(\frac1{p_0}-\frac1{\theta})})e^{-\alpha_1t}.\nonumber
\end{align*}
On the other hand, from  Lemma \ref{lemma3.18}, Lemma \ref{lemma3.20} and \eqref{3.79}, it follows that
\begin{align}
I_2&=C_4\int_0^t(1+(t-s)^{-\frac12-\frac{N}{2}(\frac{1}{q_0}-\frac{1}{\theta})})
e^{-\lambda_1(t-s)}\|m-\rho-m_\infty\|_{L^\infty(\Omega)}\|u\|_{L^{q_0}(\Omega)}ds\nonumber
\\
&\leq 2C_4M_3
M_5\varepsilon^{2}\int_0^t(1+(t-s)^{-\frac12-\frac{N}{2}(\frac{1}{q_0}-\frac{1}{\theta})})
e^{-\lambda_1(t-s)}\nonumber (1+s^{-\frac{N}{2p_0}})e^{-\alpha_1s}
(1+s^{-\frac12+\frac{N}{2q_0}})e^{-\alpha_2 s}ds\nonumber
\\
&\leq6C_4M_3
M_5\varepsilon^{2}\int_0^t(1+(t-s)^{-\frac12-\frac{N}{2}(\frac{1}{q_0}-\frac{1}{\theta})})
(1+s^{-\frac12+\frac{N}{2}(\frac{1}{q_0}-\frac{1}{p_0})})\nonumber e^{-\lambda_1(t-s)}e^{-(\alpha_1+\alpha_2) s}ds\nonumber
\\
&\leq6C_{10}C_4M_3
M_5\varepsilon^{2} e^{-\min\{\lambda_1, \alpha_1+\alpha_2\}t}(1+t^{\min\{0,\frac{N}{2}
(\frac{1}{\theta}-\frac{1}{p_0})\}})\nonumber
\\
&\leq \frac{\varepsilon }{4}(1+t^{-\frac{N}{2}(\frac{1}{p_0}-\frac{1}{\theta})})e^{-\alpha_1 t}.\nonumber
\end{align}
Combining the above inequalities, we  arrive at
$
\|(\rho-m)(\cdot,t)-e^{t\Delta}(\rho_0-m_0)\|_{L^\theta(\Omega)}\leq \frac{\varepsilon}{2}(1+t^{-\frac{N}{2}(\frac{1}{p_0}-\frac{1}{\theta}
)}) e^{-\alpha_1 t},
$ and thus complete the proof of this lemma.

By the above lemmas, we can claim that $T=T_{max}$. Indeed, if  $T<T_{max}$, by Lemma \ref{lemma3.23}, Lemma \ref{lemma3.22} and  Lemma \ref{lemma3.21}, we have
{\setlength\abovedisplayskip{4pt}
\setlength\belowdisplayskip{4pt}
$
\|(m-\rho)(\cdot,t)-e^{t\Delta}(m_0-\rho_0)\|_{L^\theta(\Omega)}\leq \frac{\varepsilon }{2}
(1+t^{-\frac{N}{2}(\frac{1}{p_0}-\frac{1}{\theta}}) e^{-\alpha_1 t},
$
$
\|\rho(\cdot,t)\|_{L^\theta(\Omega)}\leq \frac{M_1}{2}\varepsilon(1+t^{-\frac{N}{2}\left(\frac{1}{p_0}-\frac{1}{\theta}\right)})e^{-\alpha_1t}
$
as well as{\setlength\abovedisplayskip{4pt}
\setlength\belowdisplayskip{4pt}
$
\|\nabla c(\cdot,t)\|_{L^{\infty}(\Omega)}\leq \frac{M_2}{2}\varepsilon\left(1+t^{-\frac12}\right) e^{-\alpha_1 t}
$}
for all $\theta\in[q_0,\infty]$ and $t\in(0,T)$, which contradict the definition of $T$ in \eqref{3.79}. Next, the further estimates of solutions are established to ensure  $T_{max}=\infty$. \vspace{-0.5em}
\begin{lemma}\label{lemma3.24}
Under the assumptions of Proposition 3.2, there exists $M_6>0$ such that
{\setlength\abovedisplayskip{4pt}
\setlength\belowdisplayskip{4pt}\begin{align*}
\|A^\beta u(\cdot,t)\|_{L^2(\Omega)}\leq \varepsilon M_6e^{-\alpha_2 t}
\quad\hbox{for}\,\, t\in(t_0,T_{max})\,\hbox{ with}\, \,t_0=\min\{\frac{T_{max}}6,1\}.\end{align*}}
\end{lemma}

\proof For any given $\alpha_2<\lambda_1'$, we can fix $ \mu\in (\alpha_2, \lambda_1')$.  From \eqref{3.83}, it follows that
{\setlength\abovedisplayskip{4pt}
\setlength\belowdisplayskip{4pt}
\begin{align*}
\|A^\beta u(\cdot,t)\|_{L^2(\Omega)}\leq\|A^\beta e^{-tA} u_0 \|_{L^2(\Omega)}
+\int_0^t\|A^\beta e^{-(t-s)A}\mathcal{P}((\rho+m-m_\infty)\nabla\Phi)(\cdot,s)\|_{L^2(\Omega)}ds.
\end{align*}}
In the first integral, we apply  Lemma 2.3 of \cite{Cao1}, which gives
{\setlength\abovedisplayskip{4pt}
\setlength\belowdisplayskip{4pt}$$
\|A^\beta e^{-tA} u_0\|_{L^2(\Omega)}\leq
C_5|\Omega|^{\frac{N-2}{2N}} t^{-\beta}e^{-\alpha_2 t}\|u_0\|_{L^N(\Omega)}
\leq C_5 |\Omega|^{\frac{N-2}{2N}} t^{-\beta}e^{-\alpha_2 t}\varepsilon
$$}
for all $t\in(0,T)$.
Next by Lemma 2.3 of \cite{Cao1}, Lemma \ref{lemma3.18} and Lemma \ref{lemma3.22},  we have
{\setlength\abovedisplayskip{4pt}
\setlength\belowdisplayskip{4pt}
\begin{align*}
&\int_0^t\|A^\beta e^{-(t-s)A}\mathcal{P}((\rho+m-m_\infty)\nabla\Phi)(\cdot,s)\|_{L^2(\Omega)}ds
\\
\leq& C_9C_5\|\nabla\Phi\|_{L^\infty(\Omega)}|\Omega|^{\frac{q_0-2}{2q_0}}
\!\int_0^t\! e^{-\mu(t-s)}(t\!-\!s)^{-\beta} (\|m(\cdot,s)\!-\!\rho(\cdot,s)\!-\!m_\infty\|_{L^{q_0}(\Omega)}\!+\!2\|\rho(\cdot,s)\|_{L^{q_0}(\Omega)})ds
\\
\leq& M_6'\varepsilon
\int_0^t e^{-\mu(t-s)}(t-s)^{-\beta} (1+s ^{-\frac{N}{2}(\frac{1}{p_0}-\frac{1}{q_0})})e^{-\alpha_1 s}ds
\\
\leq&  M_6' \varepsilon C_{10}
  (1+t^{-1})e^{-\alpha_2 t},
\end{align*}}
where $M_6'=(M_3+M_1)C_9C_5\|\nabla\phi\|_{L^\infty(\Omega)}|\Omega|^{\frac{q_0-2}{2q_0}}$.
Therefore  there exists $M_6>0$ such that $\|A^\beta u(\cdot,t)\|_{L^2(\Omega)}\leq \varepsilon M_6 e^{-\alpha_2 t}$ for $t\in(t_0,T_{max})$.

\begin{lemma}\label{lemma3.25}
Under the assumptions of Proposition 3.2, there exists $M_7>0$ such that
$
\|c(\cdot,t)-m_\infty\|_{L^{\infty}(\Omega)}\leq M_7 e^{-\alpha_2 t}
$
for all $(t_0,T_{max})$ with $ t_0=\min\{\frac{T_{max}}6,1\}$.
\end{lemma}

\proof From \eqref{3.82} and Lemma 1.3 of \cite{Winkler7}, we have
\begin{align}
\|(c-m_\infty)(\cdot,t)\|_{L^\infty(\Omega)}
&
\leq C_1e^{-t}\|c_0-m_\infty\|_{L^\infty(\Omega)}+\int_0^t\|
e^{(t-s)(\Delta-1)}(m-m_\infty)(\cdot,s)\|_{L^\infty(\Omega)}ds\nonumber
\\
&\quad+\int_0^t\|e^{(t-s)(\Delta-1)}u\cdot\nabla c(\cdot,s)\|_{L^\infty(\Omega)}ds. \label{3.90}
\end{align}
By Lemma 1.3 of \cite{Winkler7}, Lemma \ref{lemma3.19},  we obtain{\setlength\abovedisplayskip{4pt}
\setlength\belowdisplayskip{4pt}
\begin{align}\label{3.90}
\int_0^t\|
e^{(t-s)(\Delta-1)}(m-m_\infty)(\cdot,s)\|_{L^\infty(\Omega)}ds\nonumber
\leq
&C_1 \int_0^t(1+(t-s)^{-\frac{N}{2q_0}})e^{-(t-s)}\|(m-m_\infty)(\cdot,s)\|_{L^{q_0}(\Omega)}ds
\\
\leq&
C_1C_{10}M_4\varepsilon e^{-\alpha_2t}.
\end{align}}
On the other hand, by Lemma 1.3 of \cite{Winkler7}, Lemma \ref{lemma3.20} and Lemma \ref{lemma3.21}, we get
{\setlength\abovedisplayskip{4pt}
\setlength\belowdisplayskip{4pt}\begin{align}\label{3.92}
\int_0^t\|e^{(t-s)(\Delta-1)}u\cdot\nabla c(\cdot,s)\|_{L^\infty(\Omega)}ds
\leq& C_1 \int_0^t(1+(t-s)^{-\frac{N}{2q_0}})e^{-(t-s)}\|u\cdot\nabla c (\cdot,s)\|_{L^{q_0}(\Omega)}ds\nonumber
\\
\leq& C_1 \int_0^t(1+(t-s)^{-\frac{N}{2q_0}})e^{-(t-s)}\|u(\cdot,s)\|_{L^{q_0}(\Omega)}\|\nabla c(\cdot,s)\|_{L^\infty(\Omega)}ds\nonumber
\\
\leq&6C_1M_5M_2C_{10}\varepsilon^2 e^{-\alpha_2t}.
\end{align}}
Therefore combining the above equalities, we arrive at the desired result.

{\bf Proof of Theorem 1.3 in the case $\mathcal{S}=0$ on $\partial\Omega$, part 2 (Proposition 3.2).}~
We now come to the final step to show that $T_{max}=\infty$.
 According to the extensibility criterion in Lemma \ref{lemma2.5},
it remains to show that there exists $C>0$ such that for $ t_0:=\min\{\frac{T_{max}}6,1\}<t<T_{max}$ {\setlength\abovedisplayskip{4pt}
\setlength\belowdisplayskip{4pt}
  $$ \|\rho(\cdot,t)\|_{L^\infty(\Omega)}+\|m(\cdot,t)\|_{L^\infty(\Omega)}
+\|c(\cdot,t)\|_{W^{1,\infty}(\Omega)}+
\|A^{\beta}u(\cdot,t)\|_{L^2(\Omega)}<C.$$}
From  Lemma \ref{lemma3.19} and Lemma \ref{lemma3.22},  there exists  $K_i>0$, $i=1,2,3$, such that
{\setlength\abovedisplayskip{4pt}
\setlength\belowdisplayskip{4pt}
\begin{align*}
\|m(\cdot,t)-m_\infty\|_{L^{\infty}(\Omega)}&\leq K_1 e^{-\alpha_1 t},
\|\rho(\cdot,t)\|_{L^\infty(\Omega)}\leq K_2e^{-\alpha_1 t},
\|\nabla c(\cdot,t)\|_{L^{\infty}(\Omega)}\leq  K_3 e^{-\alpha_1 t}
\end{align*}}
for $t\in(t_0,T_{max})$.
Furthermore,
 Lemma \ref{lemma3.25}  implies that
$\|c(\cdot,t)-m_\infty\|_{W^{1,\infty}(\Omega)}\leq K_3'e^{-\alpha_2t}$
with some $K_3'>0$ for all $t\in(t_0,T_{max})$.
Since $D(A^\beta)\hookrightarrow L^\infty(\Omega)$ with $\beta\in(\frac{N}{4},1)$, it  follows from Lemma \ref{lemma3.24} that
$
\|u(\cdot,t)\|_{L^\infty(\Omega)}\leq K_4 e^{-\alpha_2 t}
$
for some $K_4>0$ for all $t\in(t_0,T_{max})$. This completes the proof of Proposition 3.2.

Before we move to the next section, we remark that the following result is also valid by suitably adjusting $\varepsilon>0$ for the larger values of $p_0$ or $q_0$.\vspace{-0.5em}
\begin{corollary}
Let $N=3$ and $\int_{\Omega}\rho_0\neq\int_{\Omega}m_0$. Further,
let $p_0\in(\frac N2, \infty)$, $q_0\in(N,\infty)$ if $\int_{\Omega}\rho_0>\int_{\Omega}m_0$,
 and
 $p_0\in(\frac {2N}3, \infty)$, $q_0\in(N,\infty)$ if $\int_{\Omega}\rho_0<\int_{\Omega}m_0$. There exists
$\varepsilon>0$ such that for any initial data $(\rho_0,m_0,c_0,u_0)$ fulfilling
\eqref{1.7} as well as
{\setlength\abovedisplayskip{4pt}
\setlength\belowdisplayskip{4pt}
$$\|\rho_0-\rho_\infty\|_{L^{p_0}(\Omega)}\leq\varepsilon,\quad  \|m_0-m_\infty\|_{L^{q_0}(\Omega)}\leq\varepsilon,
\quad\|\nabla c_0\|_{L^{N}(\Omega)}\leq\varepsilon, \quad\|u_0\|_{L^{N}(\Omega)}\leq\varepsilon,$$}
\eqref{1.1} admits a global classical solution $(\rho,m,c,u,P)$.
Moreover, 
for any $\alpha_1$
$\in(0,\min\{\lambda_1, m_\infty+\rho_\infty\})$, $\alpha_2\in(0,\min\{\alpha_1,\lambda_1',1\})$, 
there exist constants $K_i$ $i=1,2,3,4$,  such that for all $t\geq 1 $   
{\setlength\abovedisplayskip{4pt}
\setlength\belowdisplayskip{4pt}
\begin{align*}
&\|m(\cdot,t)-m_\infty\|_{L^\infty(\Omega)}\leq K_1e^{-\alpha_1 t},\quad
\|\rho(\cdot,t)-\rho_\infty\|_{L^\infty(\Omega)}\leq K_2e^{-\alpha_1 t},
\\
&\|c(\cdot,t)-m_\infty\|_{W^{1,\infty}(\Omega)}\leq K_3e^{-\alpha_2t},\quad
\|u(\cdot,t)\|_{L^\infty(\Omega)}\leq K_4 e^{-\alpha_2 t}.
\end{align*}}
\end{corollary}\vspace{-1em}
\section{Proof of main results for general $\mathcal{S}$}\vspace{-1em}

In this section, we give the proof of our results for the general matrix-valued $\mathcal{S}$.
This is accomplished by an approximation procedure. In order to make the previous results applicable,
we introduce a family of smooth functions
$\rho_\eta\in C_0^\infty(\Omega)$ and $0\leq\rho_\eta(x)\leq1$ for $\eta\in(0,1),$ $\lim_{\eta\to0}\rho_\eta(x)=1$
and 
let
$\mathcal{S}_\eta(x,\rho,c)=\rho_\eta(x)\mathcal{S}(x,\rho,c).$
Using this definition, we  regularize  \eqref{1.1} as follows
 \begin{equation}\label{4.1}
\left\{
\begin{array}{ll}
(\rho_\eta)_t+u_\eta\cdot\nabla\rho_\eta=\Delta\rho_\eta-\nabla\cdot(\rho_\eta \mathcal{S}_\eta(x,\rho_\eta,c_\eta)\nabla c_\eta)-\rho_\eta m_\eta,
\\
(m_\eta)_t+u_\eta\cdot\nabla m_\eta=\Delta m_\eta-\rho_\eta m_\eta,
\\
(c_\eta)_t+u_\eta\cdot\nabla c_\eta=\Delta c_\eta-c_\eta+m_\eta,
\\
(u_\eta)_t=\Delta u_\eta-\nabla P_\eta+(\rho_\eta+m_\eta)\nabla\phi,\quad\nabla\cdot u_\eta=0,\\
\displaystyle\frac{\partial \rho_\eta}{\partial\nu}=\frac{\partial m_\eta}{\partial\nu}
=\frac{\partial c_\eta}{\partial\nu}=0,~ u_\eta=0
\end{array}
\right.
\end{equation}
with the initial data
{\setlength\abovedisplayskip{4pt}
\setlength\belowdisplayskip{4pt}
\begin{align}
\rho_\eta(x,0)=\rho_0(x),~m_\eta(x,0)=m_0(x),~c(x,0)=c_0(x),~\hbox{and}~u_\eta(x,0)=u_0(x),\quad x\in\Omega.
\label{4.2}
\end{align}}
It is observed that $\mathcal{S}_\eta$ satisfies the additional condition $\mathcal{S}=0$ on $\partial\Omega$.
Therefore based on the discussion in Section 3, under the assumptions of Theorem 1.1 and Theorem 1.3, the problem (4.1)-(4.2) admits a global classical solution $(\rho_\eta,m_\eta,c_\eta,u_\eta, P_\eta)$ that satisfies 
{\setlength\abovedisplayskip{4pt}
\setlength\belowdisplayskip{4pt}
\begin{align*}
\|m_\eta(\cdot,t)-m_\infty\|_{L^\infty(\Omega)}\leq K_1e^{-\alpha_1 t},\quad
\|\rho_\eta(\cdot,t)-\rho_\infty\|_{L^\infty(\Omega)}\leq K_2e^{-\alpha_1 t},
\\
\|c_\eta(\cdot,t)-m_\infty\|_{W^{1,\infty}(\Omega)}\leq K_3e^{-\alpha_2t}, \quad
\|u_\eta(\cdot,t)\|_{L^\infty(\Omega)}\leq K_4 e^{-\alpha_2 t}.
\end{align*}}
 for some constants $K_i$, $i=1,2,3,4$, and  $t\geq 0$. Applying a standard procedure such as in Lemma 5.2 and Lemma 5.6
of \cite{Cao1}, one can  obtain a subsequence of $\{\eta_j\}_{j\in \mathbb{N}}$ with $\eta_j\to 0$ as $j\to \infty$ such that
$
\rho_{\eta_j}\rightarrow \rho, ~m_{\eta_j}\rightarrow m, ~c_{\eta_j}\rightarrow c, u_{\eta_j}\rightarrow u \quad
\hbox{in}~ C_{loc}^{\vartheta,\frac{\vartheta}2}(\overline\Omega\times (0,\infty))
$
 as $j\rightarrow \infty$ for some $\vartheta\in (0,1)$.
Moreover, by the arguments as in Lemma 5.7, Lemma 5.8 of \cite{Cao1}, one  can also  show that $(\rho,m,c,u, P)$ is a classical solution of (1.1)
with the decay properties asserted in Theorem 1.2 and Theorem 1.3. The proofs of Theorems 1.1--1.3 are thus complete.
\vspace{-1em}
\section{Acknowledgments}\vspace{-1em}
The author Jing Li is grateful to Chinese University Hong Kong for its hospitality in January 2018, when this work was initiated.
This work is partially supported by the NSFC grants 11571363 and 61620106002, and the NUS AcRF grant R-146-000-249-114.

\end{document}